\documentclass{article}
\usepackage{amsfonts,amsmath,amscd,amssymb,amsthm}
\usepackage{stackrel}
\usepackage{fullpage}
\usepackage{xcolor}
\usepackage[T1]{fontenc}
\usepackage{imakeidx}
\usepackage{subcaption}
\makeindex[columns=3 title=Alphabetical Index]
\makeindex
\usepackage[T2A,T1]{fontenc}
\usepackage[utf8]{inputenc}
\usepackage[all]{xy}

\usepackage[alphabetic,initials]{amsrefs}
\input xy%
\usepackage{fullpage}
\usepackage{amsrefs}
\makeindex

\def\convf{\hbox{\space \raise-2mm\hbox{$\textstyle      \bigotimes \atop \scriptstyle \omega$} \space}}
 
\def\0{{\bar 0}}
 
\def\1{{\bar 1}}
\def\a{{\mathtt a}}

\def\c{{\mathtt c}}

\def\q{{\mathtt q}}

\def\Z{{\mathbb Z}}

\def\N{{\mathbb N}}

\def\R{{\mathbb R}}

\def\X{{\mathbb X}}

\def\Div{{\operatorname{Div_\Z}}}
\def\int{{\operatorname{Int_\Z}}}

\def\ds{\stackrel{\cdot}{\cup}}

\def\tors{{\operatorname{tors}}}

\def\Qbb{{\operatorname{\frak{Qco}\; \X_{\gb,\gb}}}}
\def\Qcc{{\operatorname{\frak{Qco}\; \X_{\gc,\gc}}}}

\def\Qcb{{\operatorname{\frak{Qco}\; \X_{\gc,\gb}}}}

\def\Obb{{\operatorname{\cO_{\gb,\gb}}}}
\def\Occ{{\operatorname{\cO_{\gc,\gc}}}}

\def\Ocb{{\operatorname{\cO_{\gc,\gb}}}}
\def\Spec{{\operatorname{Spec}\;}}

\def\ad{{\operatorname{ad \;}}}
\def\ado{{\operatorname{\overline{ad} \;}}}

\def\hgt{{\operatorname{ht}}}

\def\inr{{\operatorname{ind}}}
\def\rer{{\operatorname{res}}}

\def\ch{{\operatorname{ch}\:}}

\def\pwg1{{\operatorname{PWG}}}
\def\pwg{{\operatorname{pwg}}}
\def\pro{{\operatorname{Proj\,}}}

\def\wt{{\operatorname{wt}}}

\def\span{{\operatorname{span}}}

\def\Hom {{\operatorname{Hom}}}
\def\Ker {{\operatorname{Ker}\;}}

\newcommand{\ttk}{\mathtt{k}}

\newcommand{\gL}{\Lambda}

\newcommand{\itema}{\item[{{\rm(a)}}]}
\newcommand{\itemb}{\item[{{\rm(b)}}]}
\newcommand{\itemc}{\item[{{\rm(c)}}]}
\newcommand{\itemd}{\item[{{\rm(d)}}]}

\newcommand{\noi}{\noindent}
\newcommand{\ga}{\alpha}
\newcommand{\gb}{\beta}
\newcommand{\gc}{\gamma}

\newcommand{\Gd}{\Delta}

\newcommand{\gd}{\delta}

\newcommand{\gs}{\sigma}

\newcommand{\go}{\omega}

\newcommand{\gt}{\tau}
\newcommand{\gz}{\zeta}

\newcommand{\gl}{\lambda}

\newcommand{\gr}{\rho}
\newcommand{\gk}{\kappa}
\newcommand{\gep}{\epsilon}
\newcommand{\gth}{\theta}
\newcommand{\op}{\oplus}

\newcommand{\A}{\mathbb A}

\newcommand{\ot}{\otimes}

\newcommand{\fg}{\mathfrak{g}}\newcommand{\fgl}{\mathfrak{gl}}
\newcommand{\fsl}{\mathfrak{sl}}\newcommand{\fpsl}{\mathfrak{psl}}\newcommand{\osp}{\mathfrak{osp}}

\newcommand{\fr}{\mathfrak{r}}

\newcommand{\fh}{\mathfrak{h}}

\newcommand{\fb}{\mathfrak{b}}

\newcommand{\fn}{\mathfrak{n}}

\newcommand{\fsp}{\mathfrak{sp}}

\newcommand{\fl}{\mathfrak{l}}

\newcommand{\ff}{\footnote}
\newfont{\eufm}{eufm10 scaled\magstep1}
\newcommand{\pd}{\partial}
\newcommand{\ci}{\circ} \newcommand{\ti}{\times}

\newcommand{\cO}{\mathcal{O}}
\newcommand{\cC}{\mathcal{C}}

\newcommand{\cB}{\mathcal{B}}

\newcommand{\cA}{\mathcal{A}}

\newcommand{\cF}{\mathcal{F}}

\newcommand{\cG}{\mathcal{G}}

\newcommand{\cR}{\mathcal{R}}

\newcommand{\ey}{\end{eqnarray}}
\newcommand{\by}{\begin{eqnarray}}
\newcommand{\nn}{\nonumber}

\newcommand{\bco}{\begin{conjecture}}
\newcommand{\ba}{\begin{alg}}
\newcommand{\ea}{\end{alg}}
\newcommand{\eco}{\end{conjecture}}
\newcommand{\bpf}{\begin{proof}}
\newcommand{\epf}{\end{proof}}
\newcommand{\bt}{\begin{theorem}}

\newcommand{\et}{\end{theorem}}

\newcommand{\br}{\begin{rem}}
\newcommand{\er}{\end{rem}}
\newcommand{\brs}{\begin{rems}}
\newcommand{\ers}{\end{rems}}
\newcommand{\bi}{\begin{itemize}}
\newcommand{\ei}{\end{itemize}}
\newcommand{\bl}{\begin{lemma}}
\newcommand{\bsul}{\begin{sublemma}}
\newcommand{\esul}{\end{sublemma}}
\newcommand{\bp}{\begin{proposition}}
\newcommand{\be}{\begin{equation}}
\newcommand{\bc}{\begin{corollary}}
\newcommand{\bexs}{\begin{examples}}
\newcommand{\eexs}{\end{examples}}
\newcommand{\bexa}{\begin{example}}
\newcommand{\eexa}{\end{example}}
\newcommand{\bex}{\begin{exercise}}
\newcommand{\eex}{\end{exercise}}
\newcommand{\btab}{\begin{tab}}
\newcommand{\etab}{\end{tab}}
\newcommand{\el}{\end{lemma}}
\newcommand{\ep}{\end{proposition}}
\newcommand{\ee}{\end{equation}}
\newcommand{\ec}{\end{corollary}}
\newcommand{\Bc}{\begin{center}}
\newcommand{\Ec}{\end{center}}
\newcommand{\bh}{\begin{hyp}}
\newcommand{\eh}{\end{hyp}}
\newcommand{\bhs}{\begin{hyps}}
\newcommand{\ehs}{\end{hyps}}
\newcommand{\bd}{\begin{dfn}}
\newcommand{\ed}{\end{dfn}}

\newcommand{\bcu}{\bigcup}
\newcommand{\bop}{\bigoplus}

\begin{document}
\newtheorem{thm}{Theorem}[section]
\newtheorem{hyp}[thm]{Hypothesis}
 \newtheorem{hyps}[thm]{Hypotheses}
  \newtheorem{rems}[thm]{Remarks}

\newtheorem{conjecture}[thm]{Conjecture}
\newtheorem{theorem}[thm]{Theorem}
\newtheorem{theorem a}[thm]{Theorem A}
\newtheorem{example}[thm]{Example}
\newtheorem{examples}[thm]{Examples}
\newtheorem{corollary}[thm]{Corollary}
\newtheorem{rem}[thm]{Remark}
\newtheorem{lemma}[thm]{Lemma}
\newtheorem{sublemma}[thm]{Sublemma}
\newtheorem{cor}[thm]{Corollary}
\newtheorem{proposition}[thm]{Proposition}
\newtheorem{exs}[thm]{Examples}
\newtheorem{ex}[thm]{Example}
\newtheorem{exercise}[thm]{Exercise}
\numberwithin{equation}{section}%
\setcounter{part}{0}
\newcommand{\drar}{\rightarrow}
\newcommand{\lra}{\longrightarrow}
\newcommand{\rra}{\longleftarrow}
\newcommand{\dra}{\Rightarrow}
\newcommand{\dla}{\Leftarrow}
\newcommand{\xra}{\xrightarrow{\sim}}
\def\sd{\operatorname{SDer}}
\def\Str{{\operatorname{Str}\;}}

\newtheorem{Thm}{Main Theorem}


\newtheorem*{thm*}{Theorem}
\newtheorem{lem}[thm]{Lemma}
\newtheorem*{lem*}{Lemma}
\newtheorem*{prop*}{Proposition}
\newtheorem*{cor*}{Corollary}
\newtheorem{dfn}[thm]{Definition}
\newtheorem*{defn*}{Definition}
\newtheorem{notadefn}[thm]{Notation and Definition}
\newtheorem*{notadefn*}{Notation and Definition}
\newtheorem{nota}[thm]{Notation}
\newtheorem*{nota*}{Notation}
\newtheorem{note}[thm]{Remark}
\newtheorem*{note*}{Remark}
\newtheorem*{notes*}{Remarks}
\newtheorem{hypo}[thm]{Hypothesis}
\newtheorem*{ex*}{Example}
\newtheorem{prob}[thm]{Problems}
\newtheorem{conj}[thm]{Conjecture}
%

\title{Geometry and coefficients of Šapovalov elements for KM  Lie superalgebras}
\author{Ian M. Musson\ff{Research partly supported by  NSA Grant H98230-12-1-0249, and Simons Foundation grant 318264.} \\Department of Mathematical Sciences\\
University of Wisconsin-Milwaukee\\ email: {\tt
musson@uwm.edu}}

\maketitle 

\begin{abstract} 
We study \v Sapovalov   elements for a symmetrizable, integrable Kac-Moody 
Lie superalgebra $\fg$. 
Let  $\gc$ be a positive  root of $\fg$ and $m$ a positive integer, with suitable conditions on the pair   $(\gc, m).$ The \v Sapovalov element
$\theta_{\gamma,m}\in U(\fb^-)	$
has the important property that if $\gl$ lies on a certain hyperplane, then
$\theta_{\gc, m} v_\lambda$ is either zero or a highest weight vector of weight $\lambda -m\gc$ in
$M(\gl)$. We give  a closed fomula for all coefficients that arise in the induction step of the 
proof. The commutative version  of this formula relates the leading terms  using maps of graded  algebras that we call  \v Sapovalov   algebras. This suggests  a  geometric setting  for the study of 
\v Sapovalov elements. 
\end{abstract}
\noi \ref{uv.1}  {Introduction}
\\ 
 \noi  \ref{sss7.1}  {Preliminaries} \\
\noi{\ref{PA.2} {Main results on coefficients}\\
\noi \ref{1s.5}} {Proofs of Theorems   \ref{1sag} and  \ref{cot}}\\
\noi \ref{g1} {Commutative version}\\ 
\noi  \ref{GF} {Geometry of Šapovalov elements} \\ 
\noi \ref{vf} Examples


\section{Introduction.} \label{uv.1} 

Using  a remarkable combination of representation theory and algebraic geometry,   \v Sapovalov   elements and the \v Sapovalov   determinant  were introduced in \cite{Sh}.  A few examples suffice to demonstrate their importance. 
The \v Sapovalov determinant is used to construct the Jantzen filtration for semisimple Lie algebras \cite{J1} and also to compute the center of the enveloping algebra of a finite dimensional Kac-Moody Lie superalgebra \cite{Gk}, \cite{Kac4}.
 \v Sapovalov elements are crucial for understanding maps between Verma modules. Fix a positive root $\gb$ and integer $m$.  Define a \index{Hyperplane ${H}_{\gb, m}$}hyperplane 
\be \label{vat}{H}_{\gb, m} = \{ \lambda \in  {\mathfrak h}^*|(\lambda+\gr, \gb) = m(\gb, \gb)/2  \}.\ee  
We postpone some of the definitions in this introduction and for simplicity assume $\fg$ is a simple Lie algebra. The Verma module $M(\gl)$ is generated by a highest weight vector $v_\gl$ of weight $\gl$.
\index{Verma module $M(\gl)$} 
A {\it \v Sapovalov element} $\theta_{\gb, m} \in U(\fb^-)$  has the property that  for all
$\gl \in H_{\gb, m}$, 
$\theta_{\gb, m} v_\lambda$ is a highest weight vector of weight $\lambda -m\gb=\mu$ in
$M(\gl)$.
\\ \\
Define a 
{\it \v Sapovalov map} between Verma modules $M(\mu)\lra M(\gl)$  by $xv_\mu \lra x\theta_{\gb, m} v_\gl $. \ff{For Lie superalgebras a \v Sapovalov map can be zero} In the  Lie algebra case any nonzero map between Vermas is a composite of \v Sapovalov maps. If $\gb_0$ is a simple root, then $\theta_{\gb_0,  m} = e_{-\gb_0}^m$. 
Thus what is of interest in the construction of \v Sapovalov  elements
is the induction step.  Suppose that the simple reflection $s_\ga$ satisfies
\be \label{vct}
\gc = s_\ga \gb = \gb +q\ga\ee for some $q>0.$ If   $\theta_{\gb,m}$ has been constructed we need to find 
$\theta_{\gc,m}$.  
\\ \\
Many of the results of this paper come in commutative and noncommutative  versions.  Let us explain what this means.  In Theorem \ref{cgw} using  two linear maps $B_U$ and $\nabla^{t}_U$ 
we give a closed formula 
for all the coefficients 
of  $\theta_{\gc,m}$ in terms of the coefficients of 
$\theta_{\gb,m}$. 
Theorem \ref{cgw} depends on a particular ordering on the set of positive roots, and thus it makes sense to pass to the commutative associated graded ring.  
\subsubsection{The symbol  map and symmetrization} \label{syb}
First a clarification. The  symmetrization map 
$\go:S(\fn^-)\lra U(\fn^-)$ is defined in Subsection \ref{ss6.4}.
For nonzero $a\in S(\fh)$, write $a=\sum_{i=0}^n a_i$ where $a_i$ is homogeneous of degree $i$ and  $a_n\neq 0.$  Then define the {\it symbol}  \index{Symbol map $\gk$} of $a$ to be $\gk(a)= a_n \in S(\fh)$.  
We use a basis 
$\{x^\ttk_{\gs}\}$ of  $S(\fn^-)$ indexed by partitions, see
\eqref{kpar}. Then any $x\in U(\fb^-)$ can be written uniquely as
\be \label{nfn}x=\sum_\gs \go(x^\ttk_{\gs})\ot a_\gs,\ee with    
$ a_\gs\in S(\fh)$. Then define $\gk(x)= \sum_\gs  x^\ttk_{\gs}\ot\gk(a_\gs)\in S(\fb^-)$.  This carries more information than the usual symbol map for $U(\fb^-)$.
\\ \\
In one of our main results, Theorem \ref{zra} we obtain 
a closed formula
for the leading terms of the coefficients 
of  $\theta_{\gc,m}$ in terms of the leading terms of the coefficients of 
$\theta_{\gb,m}$, and three algebra maps.  We give a short summary.
\subsubsection{Algebras related to $S(\fh)$} \label{ta2} 
Fix an integer $m$. 
By \eqref{vat} $s_\ga$ maps $
{H}_{\gb, m}$ to ${H}_{\gc, m}$. Then  ${s_\ga} $ induces maps $P_{\ga} :S(\fh)\lra S(\fh)$ and 
$\cO(H_{\gc, m}) \lra \cO(H_{\gb, m}) $
such that $P_{\ga}  \phi(\mu) =\phi({s_\ga}  \cdot \mu)$.  
 Also 
  $\cO(\fh^*) = S(\fh)$ and we have a commutative diagram
\Bc
\xymatrix{&&&
 &S(\fh)\ar@{>}[d]\ar@{>}[rr]^{P_{\ga} }
&
&S(\fh)\ar@{>}[d]\\
&&&&\cO(H_{\gc, m}) \ar@{>}[rr]
&&\;
\cO(H_{\gb, m})
&
}\Ec
where the downward arrows are restrictions of functions to a hyperplane. 
The map   $P_\ga$
\index{The map $P_\ga$}
is an  important ingredient in both Theorem \ref{cgw} and Theorem \ref{zra}.

\subsubsection{Algebras related to $S(\fn^-)$ and the algebras $R(\gt,\gs)$} \label{ta1} 
\noi Fix  $\gt, \gs\in \{\gb,\gc\}$.  The 
{\it $\gt$-Šapovalov algebra} Š$_\gt$ \index{$\gt$-Šapovalov algebra} is defined by taking the direct sum weight spaces over all nonnegative integers $m$, of weight $-m\gt$ in  the localization $S(\fn^-)_e$, 
see  Section \ref{ve} for more details. 
 The  
{\it small $\gt$-Šapovalov algebra}   $\mbox{Š}^0_{\gt}$
\index{$\gt$-Šapovalov algebra!small}  is the graded subalgebra of $\mbox{Š}_{\gt}$ generated by the degree 1 part $\mbox{Š}_{\gt}(1)$. 
 The analogs 
$B_S, \nabla^t_S$ of $B_U$ and $\nabla^t_U$ are isomorphisms of
 graded  algebras
\be\label{781}\mbox{Š}_\gb  \stackrel{\nabla^t_S}{\lra} \mbox{Š}_\gb \stackrel{B}{\lra} \mbox{Š}_\gc.\ee
 Here $\nabla^t_S = \exp(td_ -)$ where $d_ -$ is a locally nilpotent derivation of
$S(\fn^-)_e$  that preserves  weight spaces.  
The map $B_S$ on  $\mbox{Š}_\gb (m)$  is multiplication by $e_{-\ga}^{mq}$.  
 Note that the multiplication takes place in $S(\fn^-)_e$. A rigorous definition of $B_S$ is given in Equation \eqref{ma4}.  The maps in  \eqref{781} restrict to isomorphisms of small Šapovalov algebras. 
The $\gt$-Šapovalov algebra is the largest algebra such that Theorem  \ref{zra} holds while $\mbox{Š}^0_{\gt}$ is the smallest.
\\ \\
Next set  
 $R(\gt,\gs) =\mbox{Š}_{\gt}^0\ot \cO(H_{\gs,0}).$
 If $\gt$ is non-isotropic,  then by shifting the origin in $\fh^*$ we can make  
$R(\gt,\gs) $ into a graded algebra 
with  degree $m$ component  $\mbox{Š}^0_{\gt}(m)\ot \cO(H_{\gs,m}),$ see Lemma \ref{td5}. \index{Algebras $R(\gt,\gs)$}
\noi Combining the  above  maps we obtain algebra isomorphisms
\be \label{zr6}  B_S\ci\nabla_S^t\ot 1: R(\gb,\gb) \lra  
R(\gc,\gb)\ee 

\be \label{zra4}  1\ot P_\ga: R(\gc,\gc) \lra  
R(\gc,\gb).\ee

\bt \label{783} If $t $ is as in  
\eqref{121d}, 
then the leading terms of the coefficients of the Šapovalov elements $\theta_{\gamma,m}=\sum_{\pi }  x_{\pi}^\ttk c_\pi$ and $\theta_{\gb,m}
=\sum_{\gs }  x_{\gs}^\ttk a_\gs$ are related by
\by 
\label{178} \sum_{\pi \in {\bf P}(m\gc)}  (1\ot P_\ga)( x_{\pi}^\ttk \gk(c_\pi))
&=& \sum_{\gs \in {{\bf P}}(m\gb)} (B_S\ci\nabla_S^t\ot 1)(x_{\gs}^\ttk  \gk(a_\gs)). 
\ey 
\et\noi
Theorem \ref{783} is a special case of  Theorem \ref{zra} which  applies to more general Lie superalgebras and uses $\Z$-forms of all the algebras involved.
\\ \\Rearranging terms in \eqref{178} we  have an isomorphism, see \eqref{2W4}  
\be \label{eeh} f: =
(1\ot P_\ga^{-1})\ci ((B_S\ci\nabla_S^t)\ot 1): R(\gb,\gb)\lra  R(\gc,\gc).\ee  such that 
$f(\gk({\theta}_{\gb,m}))=\gk({\theta}_{\gc,m}).$ 
\\ \\  
This paper is organized so that the results on coefficients can be presented as soon as possible.
\noi The next Section contains some preliminaries. After stating the key Lemma \ref{11768},  we sketch the induction step in the construction of \v Sapovalov elements. We also explain how to zoom in on the induction step, which allows us to greatly simplify notation at a cost of losing information from previous steps. 
The main results concerning coefficients are stated in Section \ref{PA.2} and proved in Section \ref{1s.5}. The proofs depend on a rather subtle cancellation property which is illustrated in Subsection \ref{Type C} in the cases of $\fsp(6)$ and $\osp(2,4)$. In 
Section \ref{g1}, apart from the results described above, we use convex geometry to show that 
the algebras Š$_\gc$  and Š$_\gb$  are finitely generated 
(provided $\fg $ is finite dimensional). We also study module categories for \v Sapovalov   algebras. 
In Section \ref{GF} we interpret our results geometrically. The final Section gives some examples.
\\ \\
I would like to thank 
Kevin Coulembier,   Maria Gorelik, 
Volodymyr Mazorchuk, Dmitriy Rumynin,  Michel van den Bergh and the referee for some helpful correspondence.

\section{Preliminaries} \label{sss7.1}
 Throughout $[n]$ denotes the set of integers 
$\{1,2,\ldots,n\}$ and 
$\ttk$ is a field of characteristic zero. Let $A$ be a $\ttk$-algebra.  
A $\Z$-{\it  form} of $A$ 
 is a subring $A_\Z $  such that 
$A_\Z\ot_\Z \ttk =A$. 

\subsection{The Lie superalgebras} \label{gsi}For $\ttk$  algebraically closed, we make the assumptions below about the Lie superalgebra  $\fg$.
\bi
\itema  $\fg=\fg(A,\gt)$ 
\index{KM Lie superalgebra $\fg(A,\gt)$} is a KM\ff{Kac-Moody Lie superalgebra. These algebras are also known as contragredient Lie superalgebras, \cite{Kac1}, \cite{M101}  Chapter 5} Lie superalgebra with symmetrizable Cartan matrix $A$
\itemb $\fg$ is integrable, that is
the adjoint action of any nonisotropic (positive or negative) simple root vector is locally nilpotent.
\itemc The enveloping algebra $U(\fg)$ \index{Integral forms!$U_\Z(\fg)$} has a $\Z$-form $U_\Z(\fg)$ that is stable under the divided power adjoint action of (positive and negative) root vectors for all simple nonisotropic roots.
\ei The  
Lie superalgebras satisfying (a)-(b) have been classified, \cite{S1}. They include all finite dimension simple Lie algebras and all basic classical simple Lie superalgebras.  However in place of $\fpsl(n|n)$ we use $\fgl(n|n)$.
\\ \\
If $\fg$
 is a simple Lie algebra, then $U_\Z(\fg)$ was constructed by Kostant based on work of Chevalley, see \cite{H}   Chapters 25-27. Fioresi and Gavarini \cite{FG} showed that hypothesis (c) also holds for classical Lie superalgebras of type A-D using an analogous construction. 
For KM Lie algebras given by a GCM, hypothesis (c) holds by \cite{Mq} Chapter 7. 
 For the Lie superalgebra $\fg=D(2,1;\ga)$ with $\ga$ irrational, there is obviously no $\Z$-form of $U(\fb^-)$.
\\ \\
Let 
 $\fg(A,\gt)$ be a KM Lie superalgebra with Cartan subalgebra $\fh$, and base  of simple roots $\Pi$.
\v Sapovalov elements need not exist for Lie superalgebras see Example \ref{2.5}, so we make the following definition. Let $W_\Pi$ be the subgroups of the Weyl group $W$ generated by all reflections corresponding to roots 
nonisotropic roots in $\Pi$. We say $\gb$ is a {\it real root}\index{Real root} if it is in the $W_\Pi$ orbit of $\Pi$. This is different from the definition in \cite{S1}.\\ \\
If  (a)-(c) hold, we show that for real roots, \v Sapovalov elements and their evaluations can  be constructed inside   the algebras $U_\Z(\fb^-)$ and $U_\Z(\fn^-)$.  In cases where (a), (b) hold but (c) does not, some results still work over $\ttk$. 
\\ \\ The assumption that $\ttk$ is algebraically closed is used only to limit the algebras considered. 
Apart from Example \ref{2.5}, we work with the fixed base $\Pi$. 

\subsection{Triangular decomposition and bilinear form}
 Here we work over $\ttk.$ 
 \index{Simple roots and positive roots $\Pi$, $\Delta^{+}$}
Let $\Delta^{+}$ be the set of positive roots  containing $\Pi$,
 and 
\be \mathfrak{g} = \mathfrak{n}^- \oplus \mathfrak{h}
\oplus \mathfrak{n}^+\nn\ee \index{Triangular decomposition of $\mathfrak{g}$}the corresponding triangular decomposition  of $\fg$. We use the Borel subalgebras $\mathfrak{b} =  \mathfrak{h}
\oplus \mathfrak{n}^+$ and
$\fb^- =
\mathfrak{n}^- \oplus \mathfrak{h}$.
  Verma modules are  induced from $\mathfrak{b}$.  We do not use other Borel subalgebras. 
\\ \\
Since $A$ is symmetrizable there is a nondegenerate invariant   bilinear form on $\fg$ \index{Non-degenerate bilinear form on $\fg$} which induces   an invariant  nondegenerate invariant symmetric bilinear form $(\;,\;)$ on $\fh^*$. For all $\ga \in \fh^*$, let $h_\ga \in \fh$ be the unique element such that \be \label{hdef}(\ga,\gb) = \gb(h_\ga) \mbox{ for all }\gb \in \fh^*.\ee Then for all $\alpha
\in \Delta^+$, choose elements $e_{\pm \alpha} \in
\mathfrak{g}^{\pm \alpha}$
 such that

\be \label{1eeh} 
[e_{\alpha}, e_{-\alpha}] = h_{\alpha}.\ee
It follows that if
$v_\mu$ is a highest weight vector of weight $\mu$ then
\be \label{ebl} e_\ga e_{-\ga} v_\mu = h_\ga v_\mu =(\mu, \ga)v_\mu . \ee
\noi 

\subsection{Weyl group}
Suppose $\Pi_{\rm nonisotropic}$  is the set of nonisotropic  simple roots.
\index{Set of nonisotropic  simple roots  $\Pi_{\rm nonisotropic}$}
  \noi
 For a nonisotropic root $\ga,$ we set $\alpha^\vee = 2\alpha /
(\alpha, \alpha)$, and denote the reflection corresponding to $\alpha$ by $
s_\ga$.      Also define 
\be \label{efh} t_{\alpha} = 2h_{\alpha}/(\alpha, \alpha).\ee For all $\gb\in \fh^*$, we have 
$t_{\alpha} (\gb) = ({\alpha}^\vee,\gb).$
The  Weyl group  $W$ is the subgroup of $GL({\fh}^*)$
generated by all reflections.  For $u \in W$ set
 \be 
N(u) = \{ \alpha \in \Delta_0^+ | u \alpha < 0 \},\qquad \ell(u) = |N(u)|.\nn\ee
We need the following well-known fact.
\bl
If $w  = s_\ga u$ with  $\ell(w)>\ell(u)$ and $\ga$ is a nonisotropic simple root, then we have a disjoint union
\begin{eqnarray} \label{Nw}
N(w^{-1}) = s_\ga N(u^{-1}) \ds \{\ga\}.
\end{eqnarray}
\el\noi 
Integrability has some important consequences. First set 
 \index{Special sets of roots ${\overline{\Delta}}^+_{0}, {\overline{\Delta}}^+_{1}$}
\be 
{\overline{\Delta}}^+_{0} = \{ \alpha \in \Delta^+_{0} | \alpha
/ 2 \not\in \Delta^+_{1}\}, \quad {\overline{\Delta}}^+_{1}= \{ \alpha
\in \Delta^+_{1}|2\alpha \not\in \Delta^+_{0}\}.\nn\ee
If $\ga$ is a nonisotropic root, then $\exp e_{\ga}\exp e_{-\ga}\exp e_{\ga}$ is a linear isomorphism  sending $\fg^\nu $ to  $\fg^{s_\ga\nu} $.
Thus $\dim \fg^\nu = \dim\fg^{s_\ga\nu} $.
The following consequence is noted in \cite{S1} Lemma 4.11.
\bl \label{xyz} If $\gb$ is a real root then
\bi \itema $\dim \fg^\gb =1$\ei
Suppose $k\gb$ is a root of $\fg$ for $k$ positive, then 
\bi  
\itemb If $\gb\in
{\overline{\Delta}}^+_{0}  \cup {\overline{\Delta}}^+_{1}$ then $k=1$.
\itemc If $\gb$ is odd nonisotropic, then  $k =1$ or $2$.
\ei
\el
\bpf Since this is true for $\gb$ simple, the result follows from the preceding remarks.\epf

\subsection{Symmetrization and derivations} \label{ss6.4} 
Let 
 $\{ U_n \}$ be the standard filtration on $U = U(\fn^-)$ and  
     let $S^n(\fn^-)$ be the graded component of the symmetric algebra     $S(\fn^-)$  of degree $n$. 
There is a linear  map  $\omega_n: S^n(\fn^-) \longrightarrow U$ defined  by
\[ \omega_n(x_1x_2 \ldots x_n)   = (1/n!)\sum_{\sigma  \in \mathcal{S}_n}
\gamma(\pi(x),\sigma ) x_{\sigma(1)}x_{\sigma(2)}\ldots
x_{\sigma(n)} \] where $\gamma(\pi(x),\sigma) =\pm 1$ is as in \cite{M101} Equation $(A.2.36)$. If $\mathbb{S}^n(\fn^-)$ denotes the  image of 
 $S^n(\fn^-)$ in $U$, we have      $U_n = U_{n-1} \oplus \mathbb{S}^n(\fn^-)$.   Let  $$\go:S(\fn^-)=\bigoplus_{n\ge 0} S^n(\fn^-) \longrightarrow \bigoplus_{n\ge 0}  \mathbb{S}^n(\fn^-) =U(\fn^-)$$ 
\index{Symmetrization map $\go$}
be  the isomorphism of
vector spaces which coincides with
 $\omega_n$ on $S^n(\fn^-)$. 
\subsubsection{Derivations} \label{taa} 
 Let $A$ be an associative algebra For  $e, a \in A$
     and all $t \in \mathbb{N} $ we have,
\be \label{1cow} e^{t}a = \sum^t_{i=0} \left( \begin{array}{c}
                t \\
                i \end{array}\right) ((\ad e)^i a)e^{r - i} .\ee
\noi Here we interpret $(\ad e)a $ as $ea-ae$.  The following consequence is well-known.
  
\bc \label{1oreset} If $e$ is a locally {\rm ad}-nilpotent nonzero divisor in $A$,  then
$\{e^n|n \in \mathbb{N}\} $
is an Ore set in $A$. We write $A_e$ for the resulting Ore localization.
\ec
\noi \noi If $\ga \in \Pi_{\rm nonisotropic}$, and  $e = e_{-\ga}.$ Then $e$ satisfies the hypothesis in 
$U = U(\fn^-)$.
\\ \\If $\ga $ is a nonisotropic root, let
$\fl$ \index{Lie superalgebra $\fl$} be the subalgebra of $\fg$ generated by $e_{\pm \ga}$ and  $h_\ga =[ e_{\ga}, e_{-\ga}]$. Thus $\fl \cong \fsl(2)$ for  $\ga$ even or $\fl \cong \osp(1|2)$ if $\ga$ is odd. 
 Let $\fr$ be the ideal of $\fn^-$\index{An ideal of $\fr$ of $\fn^-$}
formed by negative weight vectors with weight different from $-\ga$. Then $\fn^-= \fr\op \ttk e_{-\ga}$.
\\ \\
There is a simple trick that allows us to treat the cases where $\ga$ is an even root and an odd nonisotropic root uniformly.
If $e$ is homogeneous in $\fg$ 
then $\pd=\ad  e$ is the  superderivation of $\fg$ defined by $(\ad  e)(x) = [e,x]$.  
Since  $U$ is  a $\Z_2$ graded algebra it has  an automorphism $s$  given by 
$$s(u_0 +u_1) = u_0 - u_1 \mbox{ for } u_i \in U_i.$$
If  $e$ is  odd, then  $\ado e=-s\ci\ad e\ci s$  is a derivation of $U$. 
Note that  $\ado e$ and $\ad e$ have the same restriction to $\fg$, compare the twisted adjoint action of Gorelik, \cite{G3}. \index{Twisted adjoint action $\ado e$}
\\ \\
Denote the extensions of $\ad e_{\pm \ga}$ or
$\ado e_{\pm \ga}$  (for $\ga$ even or odd respectively) to derivations 
$\pd_S^{\pm}
,  \pd_U^{\pm}$ of $S(\fr^-)$, $ U(\fr^-)$. 
Then extend these maps to derivations of $S(\fn^-)$  or $U(\fn^-)$ that vanish on $e_{-\ga}$. Observe that if $x\in U(\fn^-)$, then $ \pd_U x 
=ex-xe$ for any parities. We often write $\pd_S,  \pd_U$ in place of   $\pd_S^{-},  \pd_U^{-}$ since they are used far more frequently than the other two derivations.

 \bl \label{npr}
The map $\omega:S(\fn^-)
     \longrightarrow U(\fn^-)$ is $\pd$  equivariant, that is  $\go\pd_S(u) = \pd_U\go(u)$.
\el \bpf This is well-known, \cite{M101} Theorem 6.4.4.
\epf\noi   
\subsection{Partitions}
\noi 
Set $Q^+=\sum_{\gz\in \Pi} \N\gz$.
If $\eta \in Q^+$, a {\it
partition} \index{Partitions} of $\eta$ is a map
$\pi: \Delta^+  \longrightarrow
\mathbb{N} $ such that
$\pi(\upsilon) = 0$ or $1$ for all isotropic roots $\upsilon$,
 $\pi(\upsilon) = 0$ for all even roots $\upsilon$ such that $\upsilon/2$ is a root, and
\be  
\sum_{\upsilon \in {{ \Delta^+}}} \pi(\upsilon)\upsilon = \eta\nn.\ee
For $\eta \in Q^+$, we denote by $\bf{{P}(\eta)}$ the set of partitions of $\eta$. If $\ga$ is a positive isotropic root, set ${\bf{P}}_\ga(\eta)=\{ \pi\in {\bf{P}}(\eta)|\pi(\ga)=0\}$, $\bf{{p}(\eta)}=|\bf{{P}(\eta)}|$ and   $\bf{{p}_\ga(\eta)}=|\bf{{P}_\ga(\eta)}|.$ \index{Partitions!sets of ${{\bf P}}, \bf{P}_\ga$} \index{Partitions!number of $\bf{{p}(\eta)},\bf{{p}_\ga(\eta)}$}
  Fix a Chevalley basis $\{e_{-\upsilon}| \upsilon\in \Gd^+\}$, $\{h_\upsilon| \upsilon \in \Pi\}$ for $\fn^-, \fh$ respectively. 
For each positive root $\upsilon$ we have a root vector 
$e_{-\upsilon}$. Let   $x_{\upsilon}$ be another copy of 
$e_{-\upsilon}$ but inside $S_\Z(\fn^-)$. 
Fix an order on the set $ \Delta^+$, and for $\pi$ a partition,
set 
\be \label{kpar} 
e^\ttk_{-\pi} = \prod_{\upsilon \in \Delta^+} {e^{\pi (\upsilon)}_{-\upsilon}}, \quad x^\ttk_{\pi} = \prod_{\upsilon \in \Delta^+} x^{\pi (\upsilon)}_{\upsilon},\ee 
the  products being taken with respect to this order.
If $\pi \in \bf{{P}(\eta)}$  define  $|\pi| = \sum_{\upsilon \in \Delta^+} \pi(\upsilon)$ and  
$\pi! =\prod_{\upsilon \in \Delta^+} {{\pi (\upsilon)!}}.$
Then 
 $|\pi|$ is the degree of $e_{-\pi}$ in the standard filtration.  It is useful to extend this notation to arbitrary elements of 
$U(\mathfrak{n}^-)$. Thus if $x\in 
U_n(\mathfrak{n}^-) \backslash 
U_{n-1}(\mathfrak{n}^-)$ we set $|x|=n.$
Also set

\be \label{negpar} e_{-\pi} = e^\ttk_{-\pi}/\pi!, \quad x_{\pi} = x^\ttk_{-\pi}/\pi!\ee \index{Partitions!elements $e_{-\pi},  x_{\pi}$
associated to}\noi   Fix $\ga\in \Pi_{\rm nonisotropic}$.  We say the order on $ \Delta^+$ is an {\it $\ga$-order} \index{$\ga$-order} if the terms involving 
$e_{-\ga}$ and 
$x_{\ga}$
occur last in 
\eqref{negpar}. The elements $e_{- \pi}$ and  $x_{\pi} $ with $\pi \in \bf{{P}}(\eta)$ 
form a $\Z$-basis of
$U_\Z(\mathfrak{n}^-)^{- \eta}$ and $S_\Z(\mathfrak{n}^-)^{- \eta}$
\index{Integral forms!$U_\Z(\mathfrak{n}^-), S_\Z(\mathfrak{n}^-)$}
 respectively.   
\\ \\
There are two special partitions of $m\gc$. Let
$\pi^0 \in {{\bf P}}(m\gc)$ be the unique partition of $m\gc$ such
that $\pi^0(\upsilon) = 0$ if $\upsilon \in \Delta^+ \backslash \Pi .$ The partition $m\pi^{\gc}$ of $m\gc$ is given by $m\pi^{\gc}(\gc)=m,$ and $m\pi^\gc(\upsilon)= 0$ for all positive roots $\upsilon$ different from $\gc.$
The partition $m\gs^{\gb}$ of $m\gb$ is defined similarly. 
This means that $e_{-m\pi^{\gc}} =e_{-\gc}^m/m!$ 
and  $e_{-m\gs^{\gb}} =e_{-\gb}^m/m!$. 
\\ \\
If { $\pi \in {{\bf P}}(m\gc), \pi \neq
m\pi^{\gc}$ and  $|\pi|\le m$} we call $\pi$   a {\it rogue partition of }$
{m\pi^{\gc}}.$ \index{Partitions!rogue}

\subsection{Coefficients of \v Sapovalov  elements} \label{s.2}
Fix $\gr\in \fh^*$ such that if $\ga$ is a simple root, then $(\gr,\ga) =(\ga,\ga)/2$. \index{The element $\gr$} The dot action of $W$ on $\fh^*$ is defined by $w\cdot \gl = w(\gl+\gr) -\gr$. \index{The element $\gr$!  dot action}
Fix a positive root $\gc$,  a  positive   integer $m$ and  an  
order on $ \Delta^+$. If $\gc$ is isotropic, assume $m=1$, and if $\gc$ is odd nonisotropic, assume that $m$ is odd.
We say that $\gth = \theta_{\gamma,m}\in U({\mathfrak n}^{-} )^{- m\gc} \ot S(\fh)$ is a \index{\v Sapovalov element $\theta_{\gamma,m}$}
{\it \v Sapovalov element for the pair} $(\gc,m)$ if it has the form
\be \label{rat}
\theta = \sum_{\pi \in {{\bf P}}(m\gc)} e^\ttk_{-\pi} 
 c_{\pi}= \sum_{\pi \in {{\bf P}}(m\gc)} e_{-\pi} 
 \c_{\pi},\ee where $\c_{-\pi} =\pi! c_{-\pi}$
\be \label{rbt} c_{\pi} \in S({\mathfrak h}), \; c_{\pi^0} = 1,\nn\ee  and
\be 
e_{\upsilon} \theta \in U({\mathfrak g})(h_{\gc} +{\gr}(h_\gc)
-m(\gc,\gc)/2)+U({\mathfrak g}){\mathfrak n}^+ , \; \rm{ for \; all }\;\upsilon \in \Pi. \nn\ee
We call the $c_{\pi}$ and $\c_{\pi}$ in \eqref{rat} the {\it coefficients} \index{\v Sapovalov element $\theta_{\gamma,m}$!coefficients of} of $\gth$.  
 For a simple Lie algebra, the existence of such elements was shown by \v Sapovalov, \cite{Sh} Lemma 1.
\\ \\
We evaluate $\theta_{\gamma,m}$ at points $\gl\in H_{\gc, m}$. 
Abusing notation we also write \eqref{rat} in the form $$
\theta = \sum_{\pi \in {{\bf P}}(m\gc)} e^\ttk_{-\pi} 
 c_{\pi}\in U({\mathfrak n}^{-} )^{- m\gc} \ot 
\cO( H_{\gc, m}).$$ Then the functions $c_{\pi}$ are uniquely determined.  
\bexa \label{2.5}{\rm \v Sapovalov elements need not exist for Lie superalgebras.  Let $\fg=\osp(3,2)$ This algebra has two systems of simple roots up conjugacy.  Namely
$$ \Pi_1 = \{\gep-\gd, \gd \}, \quad \Pi_2 = \{\gd-\gep, \gep \}.$$
Here  $\gep$ is an even root $\gd$ is odd nonisotropic and $\pm (\gep-\gd)$ are isotropic. Suppose we try to find a \v Sapovalov element for the pair $(\gep,m)$ using $\Pi_1$. Note that $\gep$ is not a real root for  $\Pi_1$. The root vectors $
e_{\gep-\gd}$ and $e_{-\gep}$ are given by the KM construction and we define $e_{-\gd}= [e_{\gep-\gd}, e_{-\gep}]$.  We choose $e_{\gd-\gep}$ so that \eqref{1eeh} holds. Suppose $(\gl,\gep)=m-1$ and let $M(\gl)$ be a Verma module constructed using $\Pi_1$. Since $\gep$ is not a simple root of $\Pi_1$,  we use an odd reflection to pass to  $\Pi_2$, multiply by  $ e_{-\gep}^m$ and use another odd reflection to return to  $\Pi_1$. This yields  
\by e_{\gep-\gd} e_{-\gep}^m e_{\gd-\gep}v_\gl &=& ( e_{-\gep}^m + me_{-\gd} e_{-\gep}^{m-1})e_{\gd-\gep}v_\gl\nn\\
&=& ((\gl,{\gep-\gd}) e_{-\gep}^m  + me_{-\gd} e_{-\gep}^{m-1}e_{\gd-\gep})v_\gl\nn \ey This is a highest weight vector, but since the coefficient of $ e_{-\gep}^m$ is not constant, there is no  \v Sapovalov element.
}\eexa
\subsection{Evaluation and a key  Lemma} \label{eva} Let $R$ be a $\ttk$-algebra. 
If  $\gl \in {\mathfrak h}^* $ we define  {\it evaluation at}  \index{Evaluation} $\gl$ to be the
map 
\be 
\varepsilon^\gl:R\otimes S({\mathfrak h})  \longrightarrow R,\;  \quad
\sum_{i} a_i \otimes b_i \longrightarrow \sum_{i} a_i  b_i(\gl).\nn\ee 
\noi For  $\gl\in\fh^*$, set
\[A(\lambda)_{0}  =  \{ \alpha \in \overline{\Delta}^+_{0} | (\lambda +\gr,
\alpha^\vee) \in \mathbb{N} \backslash \{0\} \}, \]
\[A(\lambda)_{1} = \{ \alpha \in \Delta^+_{1} \backslash \overline{\Delta}^+_{1} |(\lambda+\gr, \alpha^\vee ) \in 2\mathbb{N} + 1 \}, \]
and  \be \label{1tar}A(\lambda) = A(\lambda)_{0} \cup  A(\lambda)_{1}.\ee
 We also need  \be \label{Bdef} B(\lambda) = \{ \alpha \in \overline{\Delta}^+_{1} | (\lambda +\gr,
\alpha) = 0 \} . \ee \index{Sets of roots $A(\lambda) ,B(\lambda)$ }

\noi We construct \v Sapovalov elements  inductively using the next Lemma. Fix an $\ga$-order on $ \Delta^+$. \index{Key lemma of \v Sapovalov} 
\begin{lemma}\label{11768}
Suppose $\gb$ is a real root, $\mu \in H_{\gb,m}, \;\alpha \in\Pi_{\rm nonisotropic} \cap A(\mu)$ and set
\be \label{121d}\gl = s_\ga\cdot\mu,\;\;\gc = s_\alpha\gb,\;\;p = (\mu+\gr, \alpha^\vee),\;\;q = (\gc, \alpha^\vee)=-(\gb, \alpha^\vee).\ee
Assume that $ q, p  \in \mathbb{N}\backslash \{0\}$, $t = mq+p$ and \begin{itemize}
\itema  $\theta' \in
U({\mathfrak n}^-)^{-m\gb}$ is such that $v = \theta'v_\mu \in
M(\mu)$ is a highest weight vector.
\itemb  If  $\ga$ is odd, 
then $q = 2$ and $p$ is odd.
\end{itemize}
\noi Then there is a unique $\theta \in U({\mathfrak n}^-)^{-m\gc}$ such that
for all $\mu \in H_{\gb,m}$
\begin{equation} \label{121nd}
e^{p + mq}_{- \alpha}\theta' = \theta e^p_{- \alpha}.
\end{equation}  
\end{lemma} 
\bpf If $\ga$ is even this  is well-known, see for example  \cite{H2} Section 4.13 or \cite{M101} Lemma 9.4.3. \ff{For the corrected proof see: https://www.ams.org/publications/authors/books/postpub/gsm-131-Errata.pdf} 
The proof  uses the representation theory of $\fsl(2)$.
If $\ga$ is odd, nonisotropic we use the representation theory of $\osp(1,2)$. We can choose elements  $x = e_{\ga},\; y= e_{-\ga}$, and  $h = [x, y]= h_{\ga}$ in $\fg$ such that  $[h,x] =x$, $[h,y] =-y$. If  $\fl=\osp(1,2)$, then $U(\fl)$ is generated by $ \{x, y, h \}.$     
By \cite{M101} Lemma 8.3.2,  $[x, y]= (x,y)h_{\ga}$, so $(x,y)=1$.  Therefore
$$1= ([h,x],y) = (h,[x,y]) = (h,h) = (\ga,\ga)$$ so $\ga^\vee =2\ga$.   
 Under the given hypothesis we may identify $M(\gl)$ with the submodule of $M(\mu)$ generated by $v_\gl=y^pv_\mu$. 
The element $v$ from (a) has weight $\mu- m\gb$. 
Since $\ad y$ is locally nilpotent there is a positive integer $\ell$ such that 
$y^\ell v \in M(\gl),$ and we may assume $\ell > p +2m$ and  $\ell$ is odd.  
We compute that
\[ x^2y^{\ell}v =  \frac{\ell -1}{4}(\ell -p-2m)y^{\ell-2}v \in M(\gl).
\] Since  $\ell > p +2m$, we deduce that $y^{\ell-2}v \in M(\gl).$ 
Repeating this argument leads to $y^{p+2m}\theta'v_\mu \in M(\gl)= U({\mathfrak n}^-)y^pv_\mu.$ This shows the existence of $\gth$. Uniqueness follows since $y^p$ is not a zero divisor in  $U({\mathfrak n}^-)$.\epf
\noi 
Set $e = e_{- \alpha}$.
If we are given $\theta_{\gb,m}$ and apply the Lemma with 
$\gth' = \theta_{\gb,m}(\mu)$  for $\mu$ in a dense subset of $H_{\gb,m}$ we obtain the \v Sapovalov element
$\theta_{\gamma,m}$ and  the relation

\be \label{12nd} 
e^t\theta_{\gb,m}(\mu) =\theta_{\gamma, m}(s_\ga\cdot\mu)e^p  %
\ee
for all $\mu\in H_{\gb,m}.$  Note that $\gb$ is isotropic iff $\gc$ is. {\it We assume from now on that $\gb$ and  $\gc$ are nonisotropic unless otherwise stated.} A real root can be  isotropic or nonisotropic.
\subsection{Deformed multiplication and zooming in}  \label{ma3}  
\noi  Set  $\gc =s_\ga \gb$ as in \eqref{vct} and $\varpi = \frac{\ga}{(\ga,\gb^\vee)} = -\frac{\ga}{(\ga,\gc^\vee)}.$ The two parts of this Subsection both involve moving the origin in $\fh^*$ through multiples of  $\varpi.$ 
\subsubsection{Deformed multiplication}   
\noi
 Fix  $\gs\in \{\gb,\gc\}$. 
The map \be \label{jjj} y_m^*:H_{\gs, 0} \lra H_{\gs, m}.\ee
 given by 
$\mu\lra \mu \pm m\varpi$ 
where the sign is $+$ (resp. $-$) if $\gs = \gb$ (resp. $\gs = \gc$)
induces a map $y_{m}:\cO(H_{\gs,m})\lra \cO(H_{\gs,0})$.
\bl \label{td5} There is an associative multiplication 
\be 
\cO(H_{\gs,m}) \ot \cO(H_{\gs,n})\lra \cO(H_{\gs,m+n}).\nn\ee
where for $\gth\in \cO(H_{\gs,m}), \phi \in \cO(H_{\gs,n})$ we have 
\be 
(\gth\cdot \phi )(y_{m+n}\mu) =\gth(y_{m} \mu)\phi (y_{n}\mu),\nn\ee 
\el \noi
\subsubsection{Zooming in}  
\noi \index{Zooming in} 
Let $I=I_{\gc, m}$ \index{Hyperplane ${H}_{\gc, m}$!ideal defining}
be the ideal of $S(\fh)$ consisting of  functions vanishing on ${H}_{\gc,m}$, and $\cO (H_{\gc, m})= S(\fh)/I.$ Consider the line 
$L =\ttk \varpi\subset \fh^*$. Since $ H_{\gc,m}$ meets $L$ at the  single point $m\varpi$,  it follows that $S(\fh) = I\op \cO(L)$. 
Thus we have a map 
$\phi:S(\fh)  \lra \cO(L)$ which factors through $\cO( H_{\gc,m})$. 
With  $h_\alpha$ as in \eqref{hdef}, define 
$$T=-\phi(2h_\ga/(\ga,\ga))+mq-1.$$  Setting $T=t$,  the combined map sends $2h_\ga/(\ga,\ga)$ to $-p-1 = (\gl,\ga^\vee) = \varepsilon^{\gl}(2h_\ga/(\ga,\ga))$ and the diagram below commutes.

\xymatrix{
&&&&\
\cO( H_{\gc,m})
 \ar@{->}[dr]_{\varepsilon^{\gl}}  \ar@{->}[rr]^{\phi}  && \ttk[T]=\cO(L)
\ar@{->}^{T\lra t
}[dl] &\\
&&&&&\ttk
} \noi
Taking  coefficients from $\ttk[T]$ lets us to zoom in on the induction step and  allows for a more conceptual approach.   Thus the map $f$ from \eqref{eeh} induces an isomorphism of graded $\ttk[T]$-algebras, with $\deg T =0$
$$\mbox{Š}_\gb^0[T] \lra \mbox{Š}_\gc^0[T].$$
However if we do this, we lose all information about previous steps in the induction.

\section{Main results on coefficients} \label{PA.2} 
 We keep the notation of Lemma \ref{11768}, but use an arbitrary order on $ \Delta^+$. 
We need the Zariski dense  subset $\gL$ of  ${H}_{{\gb_0},m}$ given by 
\index{Zariski dense  set $\gL$ }
\be 
\Lambda =
\{\nu \in {H}_{{\gb_0},m}| \Pi_{\rm nonisotropic}\subseteq A(\nu)
\}
.\nn\ee  
Any $\gl\in \Lambda$ satisfies some integrality and positivity conditions, which allow us to apply Lemma \ref{11768} repeatedly. We assume that \be 
\gc = w{\gb_0} \mbox{ for a simple root } {\gb_0} \mbox{ and } w \in W_\Pi\nn\ee
and that $\ell(w)$ is minimal. For $\alpha \in N(w^{-1}),$ we define 
\be \label{iv} q(w,\ga) = (w\gb_0, \ga^\vee).\ee 
 If $\Pi $ is the set of simple roots, and $\gc = \sum_{\upsilon \in \Pi}
 a_\upsilon\upsilon,$ then the {\it height} \index{Height of a root} $\hgt \gc$ of $\gc$ is defined to be
$\hgt \gc =\sum_{\upsilon \in \Pi}
 a_\upsilon.$
\noi
In the Theorem below we assume $\fg$ is a a Lie superalgebra 
satisfying (a)-(c) from Section \ref{gsi}.
We require $\gc$ to be a real  root such that 
\bi 
\itema
 If $\gc$ is isotropic, then $m=1.$
\itemb If $\gc$ is odd nonisotropic, then $m$ is odd and $q=2$. 
\ei  

\bt \label{1sag} For all
$\gl \in w\cdot  \Lambda$ we have
\begin{equation} 
\theta_{\gamma,m}(\gl) = \sum_{\pi \in {{\bf P}}(m\gamma)} e^\ttk_{-\pi}c_{ \pi}(\lambda) \in
U_\Z({\mathfrak n}^-)^{- m\gamma} ,\nn
\end{equation}
where 
the coefficients $c_{\pi}(\gl) \in \ttk $ depend polynomially on $\gl \in
w\cdot  \Lambda$, $c_{\pi^0} = 1$ 
\be \label{10a}
 \deg  \; c_{\pi}(\gl) \leq m\hgt \gc - |\pi|,\nn\ee
and \ $c_{m\pi^{\gamma}}(\gl)$ is a polynomial function of
$\gl$ of degree $m(\hgt \gc - 1)$ with highest term equal to $c\prod_{\gs
\in N(w^{-1})}(\gl,\gs)^{mq(w,\gs)}$ for a nonzero constant $c$.
\et 
\noi 
If $\fg$ is a simple Lie algebra this result is due to Kumar and Letzter, \cite{KL}. \\ \\
By induction we use an $\ga$-order to  write  \begin{equation} \label{1x7} \theta_{\gb,m} =
\sum_{\gs \in {{\bf P}}(m\gb)} e^\ttk_{-\gs}a_\gs=
\sum_{\gs \in {{\bf P}}(m\gb)} e_{-\gs}\a_\gs.
\end{equation} 
where $\a_\gs =\gs! a_\gs$ and the coefficients   $a_\gs$  satisfy analogs of the conditions  in Theorem \ref{1sag}.
\\ \\ We also have   the following closed fomula for all coefficients.
\bt \label{cot} With the same hypotheses as Theorem \ref{1sag}, the coefficients $c_{ \pi}(s_\ga\cdot\mu)$  of the Šapovalov element $\theta_{\gamma,m}$ as in \eqref{rat} are given by

\by \label{ceq}  \c_{ \pi}(s_\ga\cdot\mu) 
=\sum_{j,\gs} 
 \q^{(j)} _{\gs}( e_{- \pi})
\a_\gs(\mu).
\ey

\et \noi
\noi In \eqref{ceq} the $\q^{(j)} _{\gs}$ are somewhat mysterious $\Z$-valued functions, which are defined using the adjoint action of an $\fsl(2)$  subalgebra of $\fg$ if $\ga$ is even, or an $\osp(1|2)$ subalgebra if $\ga$ is odd.
\\ \\
The integral form of 
$S(\mathfrak{h}) $ \index{Integral forms!$S_\Z(\mathfrak{h})$}
has $\Z$-basis consisting of all products $\prod_{\gz\in \Pi} \left(\begin{array}{c} h_\gz \\ c_\gz \end{array}\right)$, where $C= (\ldots c_\gz \ldots) \in \N^{\Pi}.$ Set 
$U_\Z(\fb^-)= U_\Z(\fn^-)\ot_\Z S_\Z(\mathfrak{h}).$
Note that no integrality properties are involved in the definition of a  Šapovalov element.  However   by \eqref{ceq} and induction 
\bc We have  $\theta_{\gc,m} \in U_\Z(\fb^-)$.
\ec
\bc \label{lter} In Theorem \ref{1sag}  $c_{m\pi^{\gamma}}$ is the unique coefficient of highest degree in $\theta_{\gc,m}.$
\ec
\bpf This follows easily from the given degree estimates, and the statements about the leading terms.\epf

\section{Proofs of Theorems   \ref{1sag} and \ref{cot}}\label{1s.5}

\subsection{$\Z $-bases for $U_\Z(\fn^-)_e^{-\eta}$ and $S_\Z(\fn^-)_e^{-\eta}.$}
Let $U=U(\fn^-)$ and  $U_e$  be the localization of $U$ at $e$. 
Set  $${\overline{\bf P}}(\eta)=\{\pi \in \bf{{P}}(\eta)| \pi(\ga)=0\}.$$
The elements $e_{- \pi}$ and  $x_{\pi} $ with  $\pi \in \bf{\overline{P}}(\eta)$  form a $\Z$-basis of 
$U_\Z(\mathfrak{r}^-)^{-\eta}$ and $S_\Z(\mathfrak{r}^-)^{-\eta}$.
 \index{Partitions!sets of ${\overline{\bf P}}$}

\bl \label{1uebasis}  \begin{itemize} \item[{}]
\itema The set $\{e_{-\pi} e^k|(k,\pi ) \in {\Z\ti\bf{\overline{P}}}(\eta-k\ga)\}$ forms a $\Z $-basis for the weight space $U_\Z(\fn^-)_e^{-\eta}.$
 \itemb If $u = \sum_{(k,\pi ) \in{\Z\ti\bf{\overline{P}}}(\eta-k\ga)} c_{(k,\pi )} e_{-\pi} e^k \in U_\Z(\fn^-)_e^{-\eta}$ with $c_{(k,\pi )} \in \ttk ,$ then $u \in U_\Z(\fn^-)
$ if and only if $c_{(k,\pi )} \neq 0$ implies $k \ge 0$.
\end{itemize}
 \noi  
\el \bpf    
Since    
 $U_\Z(\fn^-)_e =U_\Z(\fr^-)[e^{\pm 1}]$ is a twisted Laurent polynomial ring  the result follows.
 \epf
\subsection{Proof of Theorem \ref{cot} } \label{cs} Until Subsection \ref{zzw} we work with an $\ga$-order on $ \Delta^+$.  
\noi The main result of this Subsection, Theorem \ref{cot} relates
coefficients  of the Šapovalov elements  $\theta_{\gb,m}$ and $\theta_{\gamma,m}$. 
To prove the Theorem, we need some preparation.
We can write \eqref{12nd} as \be \label{a2} e^{t}_{- \alpha}\theta_{\gb,m}(\mu)  
= \theta_{\gamma,m}(s_\ga\cdot\mu) e^{t-mq}_{- \alpha}.\ee
 For  $e=e_{- \alpha}\in \fn^-, {\gs \in {{\bf P}}(m\gb)} $
     and $t \in \mathbb{N} $, then by \eqref{1cow}
\by \label{3cow} e^{t}e_{-\gs} &=&  \sum^t_{j=0} \left( \begin{array}{c}
                t \\
                j \end{array}\right) (\pd^j e_{-\gs})e^{mq - j} e^{t-mq} .\ey
The next result will be proved in Subsection \ref{cst}.

\bl  \label{1c} For all $\gs,j$ as above there exist $
R _{\gs, j}\in U(\fn^-)_e^{-m\gc}$
 such that \bi\itema 
\by\label{a1}\left(
\begin{array}{c}
t \\
j
\end{array}
\right)  \pd_U^{j}(e_{-\gs})e^{mq-j} + R _{\gs, j} \in U_\Z(\fn^-)^{-m\gc} 
\ey
\itemb  For all   
$\mu\in H_{\gb,m}$ 
$$\sum_{j \geq 0}  
\sum_{\gs\in {\bf P}(m\gb)} R _{\gs, j}
\a _{\gs}( \mu) =0.  $$
\ei
\el
\noi {\it Proof of Theorem \ref{cot}.} By \eqref{a1} there are unique  $\Z$-valued coefficient functions $\q^{(j)} _{\gs}$
\index{Functions $\q^{(j)} _{\gs}$} such that 
\by
\label{A88}\left(
\begin{array}{c}
t \\
j
\end{array}
\right)  \pd_U^{j}(e_{-\gs})e^{mq-j}+R _{\gs, j} &=&\sum_{\pi
\in {{\bf P}}(m\gc)} 
\q^{(j)} _{\gs}( e_{- \pi})e_{- \pi}.  
\ey 
Hence for all   
$\mu\in H_{\gb,m}$ 
$$\sum_{j \geq 0,\;\gs\in {\bf P}(m\gb)} \left(
\begin{array}{c}
t \\
j
\end{array}
\right)  \pd_U^{j}(e_{-\gs})e^{mq-j}
\a _{\gs}( \mu) 
= \sum_{\pi
\in {{\bf P}}(m\gc)} 
q^{(j)} _{\gs}( e_{- \pi})e_{- \pi}\a _{\gs}( \mu). $$ 
 \noi 
However by \eqref{a2}, \eqref{1x7} and \eqref{3cow},  the left side is also equal to 
$$e^{t}_{- \alpha}\theta_{\gb,m}(\mu) e^{mq-t}_{- \alpha}  
= \theta_{\gamma,m}(s_\ga\cdot\mu).$$ Thus 
\by \label{cxq}  c_{ \pi}(s_\ga\cdot\mu) 
=\sum_{\;j \geq 0}  \sum_{\gs \in {{\bf P}}(m\gb)}
 \q^{(j)} _{\gs}( e_{- \pi})
a_{\gs}(\mu)\nn
\ey as required. $\hfill \Box$

\subsection{The cancellation step}  \label{cst}

\noi
{\it Proof of Lemma \ref{1c}.} We use representation theory, so assume $t=p+mq$  as in Lemma \ref{11768}.
Assume  as in  Equation 
\eqref{1x7} that   $$\theta_{\gb,m} =
\sum_{\gs \in {{\bf P}}(m\gb)} e_{-\gs}^\ttk a_\gs \in U(\fb^-).
$$
Fix $N$ large enough so that
 $\pd_U^{N+1} e_{-\gs}^\ttk = 0,$ for all $\gs \in {{\bf P}}(m\gb).$ 
By Equations 
\eqref{1x7} and \eqref{3cow}
\by \label{08x} e^{t} \theta_{\gb,m}(\mu) 
e^{N-t}
 & =& 
\sum_{\gs \in {{\bf P}}(m\gb)} e^{t}e_{-\gs}^\ttk a_\gs(\mu)e^{N-t}
\nn\\&=& \sum_{\;j \geq 0}  \sum_{\gs \in {{\bf P}}(m\gb)}\left(
\begin{array}{c}
t \\
j
\end{array}
\right)  \pd_U^{j}(e_{-\gs}^\ttk)e^{N-j} 
a_{\gs}(\mu). 
\ey 
Each term on the right side of \eqref{08x} is contained in $U(\fn^-)^{-m\gb -N \ga }.$
Hence we can write 
\begin{eqnarray} 
\left(
\begin{array}{c}
t \\
j
\end{array}
\right)  \pd_U^{j}(e_{-\gs}^\ttk)e^{-j}e^{N} 
=\sum_{\pi \in {\overline{\bf P}}(m\gb+N\ga)}
 q^{(j)} _{\gs}( e_{-\pi}^\ttk)e_{-\pi}^\ttk.
\nn\end{eqnarray}

\bl \label{109b}  If  $
 q^{(j)} _{\gs}( e_{-\pi}^\ttk)\neq 0,$ 
then 
\begin{eqnarray} 
\pi(\ga)\ge N- mq.\nn
\end{eqnarray}
\el
\bpf By \eqref{08x}  and \eqref{121nd}, 
\by  
\sum_{j,\gs,\pi} 
 q^{(j)} _{\gs}( e_{-\pi}^\ttk)e_{-\pi}^\ttk
a_{\gs}(\mu) &=& e^{t} \theta_{\gb,m}(\mu)
e^{-p}e^{N-mq}\nn\\ & =& \theta_{\gc,m}(s_\ga \cdot \mu)
e^{N-mq} \in Ue^{N-mq}.\nn
\ey
This implies the result.
\epf\noi 
Hence for some 
$\cR _{\gs, j}\in U(\fn^-)_e^{-m\gc}$
we have 
\begin{eqnarray} 
\left(
\begin{array}{c}
t \\
j
\end{array}
\right)  \pd_U^{j}(e_{-\gs}^\ttk)e^{-j}e^{N} 
+\cR _{\gs, j}\in Ue^{N-mq}
\nn\end{eqnarray}

\noi If follows that Lemma \ref{1c} holds with 
$R _{\gs, j}=\cR _{\gs, j} e^{mq-N}/\gs!$. 
\subsection{Related results} \label{Cons}\noi
\noi The next result is used to obtain  upper bounds on the degree  of the coefficients of \v Sapovalov elements. 
\bl \label{A7}
If  $
 \q^{(j)} _{\gs}( e_{- \pi})\neq 0,$ then
\begin{eqnarray} \label{198x}
|\pi| \leq |\gs| + mq - j.\ey
\el \bpf All terms 
$\left(
\begin{array}{c}
t \\
j
\end{array}
\right)  \pd_U^{j}(e_{-\gs})e^{mq-j} $ have 
the same degree in the standard filtration on $U({\mathfrak n}^-)_e$ including those that have poles.  
So the result holds by 
\eqref{A88}.
 \epf \noi
\bl\label{chq} If $\q^{(j)} _{\gs}( e_{-m\pi^\gc}) \neq 0,$ then ${\gs}={{m\gs^\gb}}$ and  $j=mq$.
\el\bpf This holds since $e_{-m\pi^\gc} e_{-\ga}^{-mq}$ is a lowest weight vector for the adjoint action of $\fl$. The corresponding highest weight vector is $e_{-m\gs^\gb}$.
\epf \noi We also note a relation between  the functions $
 \q^{(j)}_\gs$ 
and an important  linear map   $B$. Let 
$B:U_e=U_\Z(\fn^-)_e \lra U_\Z(\fn^-)_e$ be the map whose restriction to $U^{-m\gb}_e$ is given  by $B(x)= x e^{mq}$. 
By \eqref{A88}\by
\label{A89}\left(
\begin{array}{c}
t \\
j
\end{array}
\right)  B(\pd_U^{j}(e_{-\gs})e^{-j})+R _{\gs, j} &=&\sum_{\pi
\in {{\bf P}}(m\gc)} 
\q^{(j)} _{\gs}( e_{- \pi})e_{- \pi}.  
\ey 

\subsection{The operator $\nabla^{t}_U$: noncommutative case}  \label{na1} 
\index{Maps $B_U$, $\nabla_U$: noncommutative case}
Let $\int$ \index{Integer valued poynomials $\int$} 
denote the ring of integer valued poynomials \index{Integer valued poynomials $\int$} in one variable. This  has $\Z$-basis  
$\small 
\left(
\begin{array}{c}
T\\
n
\end{array}
\right)$, \cite{H} Lemma 26.1. We 
 define a $\int$-linear map $\nabla=\nabla_U:U_e\ot_\Z \int\lra U_e\ot_\Z \int$ 
by
\by \label{def nabla}\nabla(x)&=&
\sum_{\;j \geq 0} 
\left(
\begin{array}{c}
T\\
j
\end{array}
\right)\pd_U^{j}(x)
 e^{-j} 
\ey
for $x\in U_e.$ Then 
for $t\in\ttk$ we obtain an operator on $U_e$ by evaluation $ \nabla^{t}_U = \nabla|_{T=t}$.  
Also  
\by
\label{199xx}
 B\ci \nabla^{t}_U(e_{-\gs})
&=& 
\sum_{\;j \geq 0}
\left(
\begin{array}{c}
t\\
j
\end{array}
\right)   \pd_U^{j}(e_{-\gs}) e^{mq-j}.
\ey

\subsection{Reordering PBW bases}\label{zzw} To prove the results on coefficients as soon as possible, we  quote a result out of sequence. 
\bc  \label{3W1}  For any rogue partition 
$\pi$ of $m\gc$ and any order on positive roots, the coefficient of $e_{-\pi}$ in $\theta_{\gc,m}$ is zero.
\ec \bpf This follows 
from 
Corollary \ref{3W} (c)
since if $a\in S(\fh)$ is nonzero, then $\gk(a)$ is nonzero.
\epf
\noi Without loss of generality we can work over $\ttk$.
  Consider  two orders on $\Delta^+$, and for $\pi \in
{\bf P}(m\gc)$, set $e^\ttk _{-\pi} = \prod_{\alpha \in \Delta^+} e^{\pi
(\alpha)}_{- \alpha}$, and $\overline{e}^\ttk _{-\pi} = \prod_{\alpha \in
\Delta^+} e^{\pi (\alpha)}_{- \alpha}$, the products being taken with
respect to the given orders. Then we have two bases $\{
e^\ttk _{-\pi}|\pi \in {{\bf P}}(m\gc)\}$  and $
\{\overline{e}^\ttk_{-\pi}|\pi \in {{\bf P}}(m\gc)\}$  for $U(\fn^-)^{-m\gc}$.  However  to write $\theta_{\gc,m}$ we do not need rogue partitions. So define
 $${{\bf Q}}(m\gc)=\{\pi \in {{\bf Q}}(m\gc)| \pi 
\mbox{ is not a rogue partition of } m\gc
\}.$$ 
Consider a basis 
$\cB = \{{b_\pi}|\pi \in {{\bf Q}}(m\gc)\}$  indexed by partitions.  
Fix a total order on the set ${{\bf Q}}(m\gc)$ such that if $\pi, \gs \in {{\bf Q}}(m\gc)$ and $|\pi| > |\gs|$ then $\pi$ precedes $\gs,$ and use this order on partitions to induce orders on the basis ${\cB}$.  If 
$x \in U(\fb^-)^{-m\gc}$ we can write $x$ uniquely in the form

\be 
x \; = \; \sum_{\pi \in {{\bf Q}}(m\gc)}
b_\pi h_\pi \nn.\ee 
We list various desirable conditions on the basis $\cB.$
   \bi \itema
$\deg  h_\pi \le m\hgt \gc - |\pi|$  
for all $\pi \in {{\bf Q}}(m\gc),$
\itemb As a polynomial function of
$\gl$,  
$h_{m\pi^{\gamma}}$ has leading term  $c\prod_{\gs
\in N(w^{-1})}(\gl,\gs)^{mq(w,\gs)}$ for a nonzero constant $c$.
 \ei  
We apply this to the 
sets ${\bf B_1} = \{
e^\ttk _{-\pi}|\pi \in {{\bf Q}}(m\gc)\}$  and ${\bf B_2} = \{\overline{e}^\ttk_{-\pi}|\pi \in {{\bf Q}}(m\gc)\}.$ 
\bl \label{der1}   We have $\span \;{\bf B_1} =\span \;{\bf B_2}$
and the change of basis matrix $C=(c_{\pi,\zeta})$ from the basis ${\bf B_2}$ to ${\bf B_1}$ 
is an upper triangular with all diagonal entries equal to $\pm1.$ \el %
\noindent \bpf This follows since 
the associated graded ring arising from the standard filtration on
$U=U({\mathfrak n}^-)$ is supercommutative.
\epf

\bexa \label{ex1}{\rm Suppose that $\fg = \osp(5,2)$.  Then $\fg$ has a basis of simple roots consisting of the roots $\ga_1= \epsilon_{1} - \epsilon_{2}$, $\ga_2=\epsilon_{2} - \delta $ and $\ga_3=\delta$.   Consider the partitions $\pi, \gs $ of  $\gc =
\epsilon_{1} + \epsilon_{2}$ taking nonzero values 
$ \pi(\epsilon_{1} - \gd) =\pi(\epsilon_{2} + \delta)=1$  
and $\gs( \epsilon_{1} + \epsilon_{2})=1.$  For different orders on the positive roots we have
 $e_{-\pi}^\ttk=e_{\gd-\epsilon_{1} } e_{-\epsilon_{2} - \delta}$,   
$\overline{e}^\ttk_{-\pi} = e_{-\epsilon_{2} - \delta}e_{\gd-\epsilon_{1} } $
and $e_\gs^\ttk =\overline{e}^\ttk_\gs =e_{ -\epsilon_{1} - \epsilon_{2}}$.  We have,  
$e_{-\pi}^\ttk + \overline{e}_{-\pi} = c\overline{e}_{-\gs}$ for a nonzero constant $c$, or in matrix form
$$\left[\begin{array}{c} e^\ttk_{-\pi}\\ {e^\ttk}_\gs \end{array}\right] =\left[\begin{array}{cc} -1&c\\ 0&1 \end{array}\right]
\left[\begin{array}{c}  \overline{e}^\ttk_{-\pi}\\ \overline{e}^\ttk_{-\gs}\end{array}\right].$$
}
\eexa\noi

\bl \label{derx}
For $x \in U({\mathfrak n}^-)^{-m\gc}
\otimes S({\mathfrak h}),$ and using the bases ${\bf B_1}$ and  ${\bf B_2}$ write
\be 
x \; = \; \sum_{\gz \in {{\bf Q}}(m\gc)}
\overline{e}^\ttk_{-\gz} f_\gz \; = \; \sum_{\pi \in {{\bf Q}}(m\gc)}
{e}_{-\pi}^\ttk g_\pi.\nn\ee  Then  the basis ${\bf B_1}$ satisfies conditions $(a)-(b)$ iff the basis  ${\bf B_2}$ 
satisfies the corresponding condition. 
\el
 \bpf  Assume ${\bf B_2}$ satisfies conditions $(a)-(b).$ Using the matrix $C$ from Lemma \ref{der1} we can write \[{e}^\ttk_{-\pi} = 
\sum_{\zeta} c_{\pi,\zeta}\overline{e}^\ttk_{-\gz}\] 
 Then 
$$f_\zeta = \sum_{\pi \in {{\bf Q}}(m\gc)}
c_{\pi,\zeta}g_{\pi}.$$
 It follows that $f_\zeta$ is a linear combination of polynomials of
degree less than or equal to $m\hgt \gc - |\gz|.$ 
Thus ${\bf B_1}$ 
 satisfies (a). Also the partition 
 $m\pi^{\gc}$ occurs last in both orderings, and ${e}^\ttk_{-m\pi^{\gc}} =e_{-\gc}^m=\overline{e}^\ttk_{-m\pi^{\gc}}$ occur with the same coefficient. Thus ${\bf B_1}$ 
 satisfies (b).\epf

\subsection{Proof of Theorem \ref{1sag}} \label{dt}  
By Lemma \ref{derx}, it is enough to prove  Theorem \ref{1sag} when we fix an $\ga$-order on $\Gd^+$.
\noi  Write \be 
 w =
s_{\ga}u, \quad   \gb = u{\gb_0}, \quad \gamma = w{\gb_0} = s_{\alpha} \gb,
\nn\ee
 with $\alpha \in \Pi_{\rm nonisotropic}$ and $\ell(w) =\ell(u) + 1.$
Suppose $\nu \in \Lambda$ and set  
\be 
\mu =
u \cdot\nu, \quad \lambda = w \cdot\nu
= s_\ga \cdot \mu.
 \ee
 By induction, $$\deg  a_{\gs}(\mu) \le m\hgt \gb - |\gs|.$$ By Theorem   \ref{cot}, 
\be \label{sas}\mbox{deg} \; c_{\pi}(\gl)  \leq  \max\{j + \; \mbox{deg} \;
a_{\gs}(\mu)\;|\; \q^{(j)} _{\gs}( e_{- \pi})
\neq 0\}.\ee
  Now if $\q^{(j)} _{\gs}( e_{- \pi})
\neq 0$ then by (\ref{198x}) 
$j \leq |\gs| +mq -|\pi|.$  %
We eliminate $j$ by combining  this inequality with \eqref{sas}, giving the first inequality  below. The second holds by induction and the last equality by  \eqref{121d}
\begin{eqnarray} 
\mbox{deg} \;\; c_{\pi}(\gl) \nonumber &\leq&
\mbox{deg} \;\; a_{\gs}(\mu)+ |\gs| + mq - |\pi|
\\ &\le& m \hgt \gb - |\pi| +mq\nonumber\\
&=&  m\hgt \gc - |\pi|.
\nn\end{eqnarray}
\subsubsection{Leading term} \label{tdt} First note that by \eqref{iv} we have 
\be 
q(w,\ga) = (\gc, \alpha^\vee) = q\nn\ee and 
\be 
q(u,\gt) = (u\gb, \gt^\vee) =  (w\gb, s_\ga \gt^\vee)= q(w,s_\ga \gt) .\nn\ee  \noi  
It follows from Theorem \ref{cot} and  Lemma  \ref{chq} that, modulo terms of lower
degree
\be 
 c_{m\pi^{\gamma}}(\gl) \equiv \left( \begin{array}{c}
t\\
mq
\end{array}   \right)
a_{m\gs^{\gb}}(\mu) .\nn\ee
 By induction 
 $a_{m\gs^{\gb}}(\mu)$ has
highest term $c'\prod_{\gt \in N(u^{-1})}(\mu,\gt)^{mq(u,\tau)}$ as a polynomial in $\mu$, for
a nonzero constant $c'$.   
Thus  as a polynomial in $\gl=s_\ga \cdot \mu$, $a_{m\gs^{\gb}}(\mu)$ has
highest term 
\by 
c'\prod_{\gt \in N(u^{-1})}(\gl,s_\ga\gt)^{mq(u,\tau)} 
&=&c'\prod_{\gs \in s_\ga N(w^{-1}), \gs\neq \ga}(\gl,\gs)^{mq(w,\gs)},\nn\ey 
On the other hand $t=\pm (\gl,\ga) + C$, with some constant $C$ independent of $\mu, \gl$.  Thus
\be 
\left( \begin{array}{c}
t\\
mq
\end{array}   \right)\equiv t^{mq}\equiv (\gl,\ga)^{mq}
\nn \ee
modulo terms of lower degree and up to a nonzero constant multiple. The
claim about the leading term of $c_{m\pi^{\gamma}}(\gl) $ in Theorem \ref{1sag} follows from the preceding
equations and \eqref{Nw}.

\subsection{Formula for coefficients}\label{cfw}
Here we work with an $\ga$-order on positive roots. 

\bt \label{cgw} If $t = mq+p$, as in  \eqref{121d}, 
then the  Šapovalov elements $\theta_{\gamma,m}$ and $\theta_{\gb,m}$ are related after evaluation on $\mu\in{H}_{\gb, m}$
by 
\by 
\label{177} \sum_{\pi \in {{\bf P}}(m\gamma)}  (1\ot P_\ga) e_{-\pi} c_\pi
&=& \sum_{\gs \in {{\bf P}}(m\gb)} ((B_U\ci\nabla_U^t)\ot 1)e_{-\gs} \a_\gs. 
\ey
\et
\bpf 
 Multiply the left side of \eqref{ceq}  by $e_{-\pi}$ and sum over 
$\pi \in {\bf P}(m\gamma)$ to obtain the first equality
Then use \eqref{A89}  and \eqref{199xx} 

\by \label{czq} \sum_{\pi \in {{\bf P}}(m\gamma)}  e_{- \pi}c_{ \pi}(s_\ga\cdot\mu) 
&=&\sum_{j, \gs, \pi} 
 \q^{(j)} _{\gs}( e_{- \pi})e_{-\pi}
\a_\gs(\mu).\nn\\
&=&\sum_{j, \gs, \pi}  \left(
\begin{array}{c}
t\\
j
\end{array}
\right) 
\pd_U^{j}(e_{-\gs})e^{mq-j}\a_\gs(\mu)\nn\\
&=& \sum_{\gs \in {{\bf P}}(m\gb)} ((B_U\ci\nabla_U^t)\ot 1)e_{-\gs} \a_\gs.\nn 
\ey This gives the result since 
$P_\ga  c_\pi(\mu) = c_{ \pi}(s_\ga\cdot\mu)$.
\epf
\noi 
We can also write 
\eqref{177} in the form 

\by 
\label{17x} \sum_{\pi \in {{\bf P}}(m\gamma)}  (1\ot P_\ga) \go(x_{\pi}) \c'_\pi
&=& \sum_{\gs \in {{\bf P}}(m\gb)} ((B_U\ci\nabla_U^t)\ot 1)\go(x_{\gs}) \a'_\gs. 
\ey where $\gk(\c'_\pi)= \gk(\c_\pi)$ and $\gk(\a'_\gs) =\gk(\a_\gs) $,
\subsection{Uniqueness of \v Sapovalov elements} \label{1s.1}

 The \v Sapovalov determinant $\det F_{\eta}$, the Jantzen filtration $\{M_{i}(\lambda) \}$ and sum formula are defined as usual, \cite{M101} Chapter 10.  
Except for a brief remark, proofs are omitted as they are 
the same as 
in \cite{Gk},   \cite{KK} or  \cite{M101}.   Recall the definitions of
$A(\lambda), B(\lambda)$ from equations (\ref{1tar}), (\ref{Bdef}).
\begin{theorem} \label{shapdet}
Up to a nonzero constant factor, $\det F_{\eta} = D_{1}D_{2}D_{3},$
where
\begin{eqnarray*}
D_{1} & = & \prod_{\alpha \in \overline{\Delta}^+_{0}} \;\;
\prod^\infty_{r=1}
            (h_{\alpha} + (\rho, \alpha) - r(\alpha,\alpha)/2)^
            {{\bf p}(\eta - r \alpha)}, \\
D_{2} & = & \prod_{\alpha \in
\Delta^+_{1}\backslash\overline{\Delta}^+_{1}}\;\;
            \prod^{\infty}_{\stackrel{\scriptstyle{r=1}}{\scriptstyle{r \; {\rm odd}}} }
            (h_{\alpha} + (\rho, \alpha) - r(\alpha,\alpha)/2)^
            {{\bf p}(\eta - r \alpha)}, \\
D_{3} & = & \prod_{\alpha \in \overline{\Delta}^{+}_{1}} \;\;
(h_{\alpha} +
             (\rho, \alpha))^{{\bf p}_{\alpha}(\eta - \alpha)}.
\end{eqnarray*}
\end{theorem}
\noi This is proved by first finding the leading term and then using representation theory to analyze $M(\gl)$ in the most generic cases.
Suppose  $\alpha$ is a positive nonisotropic root of $\mathfrak{g}$, and $r \in
\mathbb{Q} ^+$ and $ (\lambda + \rho, \alpha^{\vee} )  =  r $.
Then either $M(\lambda)$ is simple, or
$M(\lambda)$ has a submodule $U$ which is isomorphic to
a finite direct sum of modules isomorphic to
$M(s_\alpha \cdot \lambda)$ and $M(\lambda)/U$
is simple, \cite{KK} Lemma 3.3 (b), 
\cite{M101} Lemma 10.2.4.  However using Lemma \ref{xyz} it follows that if $M(\lambda)$ not simple, then $U\cong M(s_\alpha \cdot \lambda)$. Later in the proof it is shown that 
$\ga\in A(\lambda)$.
\\ \\
We now state the Jantzen sum formula for the Verma
modules $M(\lambda)$.  
\begin{theorem} \label{Jansum}
For all $\lambda \in \mathfrak{h}^*$
\[ \sum_{i > 0} \ch M_{i}(\lambda) = \sum_{\alpha \in A(\lambda)}
\ch M(s_\alpha \cdot  \lambda) + \sum_{\alpha \in
B(\lambda)} \epsilon^{\lambda -\alpha}p_\ga.\]
\end{theorem}
\bt \label{17612} Suppose $\gth_1, \gth_2$ are \v Sapovalov elements for the pair $(\gc,m)$.  Set $H = H_{\gc, m}.$ Then
for all $\gl \in H$ we have $\gth_1 v_\gl =\gth_2 v_\gl$. 
\et \bpf
Let $\{{M}_{i}(\lambda)\}_{i\ge 1}$ be the Jantzen filtration on 
${M}(\lambda)$. Then ${M}_{1}(\lambda)$ is the unique maximal submodule of ${M}(\lambda)$.  
Set
\[\gL = \{\gl \in H| A(\gl) = \{\gc\}, \;B(\gl) = \emptyset \},\]
if $\gc$ is nonisotropic, and
\[\gL = \{\gl \in H| B(\gl) = \{\gc\}, \;A(\gl) = \emptyset \},\]
if $\gc$ is isotropic.  If $\gl \in \gL$
it follows from Theorem \ref{Jansum}, that
$ M_1(\gl)^{\gl-m\gc}$ is one-dimensional.
Because ${M}_{1}(\lambda)$ is the unique maximal submodule of
${M}(\lambda)$, $\gth_1 v_\gl$ and $\gth_2 v_\gl$ are proportional. Then from the requirement that $e_{-\pi^0}$ occurs with coefficient 1 in
a \v{S}apovalov element we have $\gth_1 v_\gl =\gth_2 v_\gl$.  Since $\gL$ is Zariski dense in $H$, the result 
follows.
\epf
\noi 
A simple example shows that uniqueness of \v Sapovalov elements fails if $\fg(A,\gt)$ is replaced by a subfactor.  Consider the matrix
$$A= \left[\begin{array}{ccc} 2&-1&0\\-1& 0&1\\0&1&-2 \end{array}\right]$$
If $\gt =\{2\}$, then $\fg(A,\gt)\cong \fgl(2|2). $ The subfactor  
$\fg=\fpsl(2|2)$ has simple roots $\ga, \gb, \gc$, where $\ga+ 2\gb+ \gc=0.$ 
If we use the Borel subalgebra with simple root vectors $e_{13}, e_{34}, e_{42},$ then $e_{13}$ and $e_{42}$ have the same weight $\ga+\gb $, so if $(\gl,\ga+\gb)=0$ we have $$\dim \Hom_{U(\fg)}(M(\gl-(\ga+\gb)), M(\gl))=2.$$ 
\noi 

\section{Commutative version} \label{g1}

\subsection{The operator $\nabla_S$: commutative case}  \label{na2}
Recall the derivations $ \pd_S^{\pm}$ of  $S(\fn^-)$ defined in Subsection \ref{taa} and that $x_{\ga}$ is a copy of $e_{-\ga}$ in $S(\fn^-)$.
Define derivations  $d_{\pm}$ of  $S(\fn^-)_e$  by 
\be\label{A2} d_{\pm} = x_{\ga}^{\pm 1} \pd_S^{\pm}.\nn\ee 
\index{Derivations $d_\pm$} 
Let $\Div$ be the ring of divided powers of $T$. \index{Ring  of divided powers $\Div$}
This has $\Z$-basis  
 $T^{(n)}= \frac{T^n}{n!}$
  for $n\ge 0$.
Define a $\Div $-linear automorphism  $\nabla_S:S(\fn^-)_e\ot \Div \lra S(\fn^-)_e\ot \Div$ by
\by \label{Sdef nabla}\nabla_S(x) = 
(\exp T d_-)(x).
\ey
for $x\in S(\fn^-)_e.$ Then 
for $t\in\ttk$ set  $ \nabla^{t}_S = \nabla_S|_{T=t}=\exp t d_-$.    
\index{Maps $B_S$ and $\nabla_S$: commutative case}  It follows from (c) in Subsection \ref{gsi}, that $\exp   t d_\pm$ preserves $S_\Z(\fn^-)_e$.
\bl\label{dra}\bi \itema $\nabla_S$ and $\nabla^t_S$ are graded algebra maps.
\itemb We have $\nabla^{t}_S\ci \nabla^{s}_S = \nabla^{t+s}_S$ for all $s,t\in \ttk$.
\itemc $( \nabla_U\ci\go)(u) \equiv (\go \ci\nabla_S)(u)$ modulo lower degrees in $T$.
\itemd For $x$ as in \eqref{nfn} we have $( \gk\ci\nabla_U\ci\go(x)) \equiv (\gk\ci\go \ci\nabla_S)(x)$ mod lower degrees in $T$.
\ei
\el \bpf  First (a) and (b) hold for the exponential of any locally nilpotent derivation.  
Since $\go\pd^j_S(u) = \pd^j_U\go(u)$ by Lemma \ref{npr} 
\by \label{dea}\nabla_U\ci
\go(u)&=&
\sum_{\;j \geq 0}  \left(
\begin{array}{c}
T\\
j
\end{array}
\right)
\pd_U^{j } (\go(u)
 e^{-j})
\nn\\ 
&=& \sum_{\;j \geq 0} \left(
\begin{array}{c}
T\\
j
\end{array}
\right)
\go(\pd_S^{j}(u))
 e^{-j}.\nn
\ey
 A similar calculation shows that 
\by \label{bea} \go\ci\nabla_S (u)
&=& \go\ci\sum_{\;j \geq 0} \frac{T^j}{j!}\pd_S^{j}(u)
 e^{-j} \nn\\
&=&
 \sum_{\;j \geq 0} \frac{T^j}{j!}\go\pd_S^{j}(u)
 e^{-j}.
\ey 
This proves (c) and  (d)  is an immediate consequence.
\epf \noi


\subsection{ \v Sapovalov subalgebras}  \label{ve} 
 Note that  $S(\fn^-)_e$  has a weight space decomposition.  Fix  $\gt\in \{\gb,\gc\}$. 
 The $\gt$-Šapovalov algebra is the subalgebra of $S(\fn^-)_e$ given by  
$$
\mbox{Š}_\gt
= \bigoplus_{m\ge 0}
\mbox{Š}_\gt(m), \quad 
\mbox{ where }\mbox{Š}_\gt(m)= (S(\fn^-)_e)^{-m\gt}.$$ 
 Note that Š$_\gt$ is a graded algebra, but if $\gt$ is not a simple root, the grading is not compatible with the usual grading on $S(\fn^-)$.   What makes up for this is that  all the algebras are graded by the root lattice and this grading behaves well. Let $\mbox{Š}^0_{\gt}$ be the subalgebra of $\mbox{Š}_{\gt}$ generated by $\mbox{Š}_{\gt}(1)$. \index{$\gt$-Šapovalov algebra!integral form} The integral forms of  $\mbox{Š}_{\gt}$ and  $\mbox{Š}^0_{\gt}$
are defined by $\mbox{Š}_{\gt,\Z} = \mbox{Š}_{\gt}\cap S_\Z(\fr)[e^{\pm 1}]$
and  $\mbox{Š}^0_{\gt,\Z}= \mbox{Š}^0_{\gt}\cap S_\Z(\fr)[e^{\pm 1}]$
respectively. 
 Both derivations  $d_\pm$
preserve weight spaces and the integral forms $\mbox{Š}_{\gt,\Z}$. 
\\ \\
From now on we restrict  $\nabla_S=\exp T d_-$ to Š$_\gb[T]$.  Obviously $\nabla_S$ is  an algebra automorphism  
of Š$_\gb[T]$.  Also there is  an isomorphism  of graded algebras 
\be\label{BSD} B_S:\mbox{Š}_\gb[T]\lra \mbox{Š}_\gc[T]\ee which is given
on Š$_\gb(m)$ by $x\lra xe^{mq}$. 


\noi \subsection{Module categories} \label{PA.101}

\subsubsection{Induction and restriction} \label{td2} 
Let $f¨:R\lra A$ be a homomorphism of rings. Denote their module categories by $R$-mod and $A$-mod respectively. For $N\in R$-mod denote the induced $A$ module $\inr_R^A\; N$. 
If $M\in A$-mod we obtain an $R$-module   $\rer_R^A \;M$ and we  have a natural isomorphism 
\be \label{GF11}\Hom_A(\inr_R^A\; N, M) 
\cong  \Hom_R(N, \rer_R^A \;M).\ee
We denote the functors $\inr_R^A\; ,\rer_R^A $ respectively by $f^*$ and $f_*$. 
Taking $M=\inr_R^A\; N$ and the identity map on the left side of \eqref{GF11},  we obtain an adjunction isomorphism of $R$-modules  
$N\lra f_* f^* N.$ Explicitly $N \lra 1\ot N  \subseteq A\ot_R N.$ Taking $R=N$ we have $R \lra 1\ot R  \subseteq A\ot_R R \cong A.$
If 
 $R=\bop_{n\ge 0} R(n)$ is  a graded $R$ ring, let ${\bf gr} R$ be the category of graded  $R$-modules.  
\subsubsection{Definition of $B_S$} \label{td1} 
We now apply this to the rings
$R=\mbox{Š}_{\gb,\Z}, S= \mbox{Š}_{\gc,\Z}$, $A=S_\Z(\fn^-)_e$ and the  inclusions $i:R\lra A$ and $j:S \lra  A$.
By adjunction we have $R \lra 1\ot\oplus_{m\ge 0}R(m)$. 
Since $i^*R$ is an $A$-module it is legitimate to form the sum
$\oplus_{m\ge 0}e^{mq}\ot R(m)$
which maps to $\oplus_{m\ge 0}e^{mq}
 R(m)
\subseteq (S(\fn^-)_e)^{-m\gc}=$Š$_\gc(m)$. It is easy to see the inclusion is an equality.  However we have suppressed the restriction  $j_*$ (as often happens).  To do this more formally, let $N= \op N(k) $ be a graded $R$-module and define a map $B_k: 1\ot N(k) \lra e^{kq}\ot N \lra e^{kq} N$.  Then $B^*: {\bf gr} R\lra {\bf gr} S$ is the map  
which restrcited to degree $k$ is given by 
\be \label{ma4}
N(k)\stackrel{i_* i^*}{\lra} 1\ot N(k)\stackrel{B_k}{\lra} e^{kq} \ot N(k) \stackrel{j_* }{\lra}  e^{kq}  N(k).\ee
In particular, taking $N=R$ gives the map $B_S: R\lra S$.
\\ \\
The functor  $(\nabla^t_S)^*$ is much easier to describe. Indeed $(\nabla^t_S)^*N$ is just the $R$ module 
$N$ twisted by the automorphism $\nabla^t_S$. Thus $(r,n)\lra \nabla^t_S(r)n$.
\subsubsection{Quotient categories} 
For a graded ring  $R=\bop_{n\ge  0} R(n)$,  the irrelevant ideal is  $R_+=\bop_{n> 0} R(n).$  Let 
 $N=\bop_{n\ge n_0} N(n)$ be a graded $R$-module. We say that $N$ is bounded if $N(n)=0$ for $n\gg 0$ and {\it torsion} if it is a direct limit of bounded modules. Let $\tors$ be the full subcategory 
${\bf gr} R$ consisting of torsion modules. 
\index{Category of graded $R$-modules ${\bf gr} R$}
\index{Category of graded $R$-modules ${\bf gr} R$!mod torsion ${\bf gr} R/\tors$}
We form the quotient category ${\bf gr} R/\tors$. 
\bt \label{HF9} If $R=\mbox{Š}^0_{\gb,\Z}, S= \mbox{Š}^0_{\gc,\Z}$, the maps $\nabla^t_S$ and $ B_S$ induce equivalences of categories 

\be \label{LF9}{\bf gr} R  \stackrel{(\nabla^t_S)^*}{\lra} {\bf gr} R  \stackrel{B^*}{\lra}{\bf gr} S\ee  
and 
\be \label{LF10}{\bf gr} R/\tors  \stackrel{(\nabla^t_S)^*}{\lra} {\bf gr} R/\tors  \stackrel{B^*}{\lra}{\bf gr} S/\tors.\ee
\et
\noi 
If $\gb$ is  isotropic then $\mbox{Š}^0_\gb$ is finite dimensional and local. Thus the irrelevant  ideal is torsion and ${\bf gr} \;\mbox{Š}^0_\gb/\tors  ={\bf gr} \;\mbox{Š}^0_\gb(0)/\tors  $.
\subsection{Proof of Theorem \ref{783}}  \label{GF2}

If $I_m=S_\Z (\fh) \cap I( H_{\gb,m})$,
 we define $\cO_\Z(H_{\gb,m})=S_\Z (\fh) /I_m$. \index{Integral forms!$\cO_\Z(H_{\gb,m}), \cO_\Z(H_{\gc,m})$}
This  is a $\Z$-form of $\cO(H_{\gb,m})$. 
For  $\gt, \gs\in \{\gb,\gc\}$ we define a \index{Algebras $R(\gt,\gs)$!
$R_\Z(\gt,\gs)$} graded ring
$R_\Z(\gt,\gs)= \mbox{Š}^0_{\gt,\Z}\ot_\Z \cO_\Z(H_{\gs,0})$
 which has for all $m\ge 0$, degree $m$ part 
$$ R_\Z(\gt,\gs)(m)= \mbox{Š}^0_{\gt,\Z}(m)\ot_\Z \cO_\Z(H_{\gs,m}),$$
where we use Lemma \ref{td5} to multiply the second factor. 
We have  homomorphisms  
\be \label{zra1} 
1\ot P_\ga:  R_\Z(\gc,\gc)\lra  R_\Z(\gc,\gb)\ee and 
\be \label{zra2}  B_S\ci\nabla_S^t\ot 1: R_\Z(\gb,\gb) \lra  R_\Z(\gc,\gb).\ee

\bt \label{zra} Assume $\gb$ is a real root. If $t = mq+p$, as in  
\eqref{121d}, 
then the leading terms of the coefficients of the Šapovalov elements $\theta_{\gamma,m}=\sum_{\pi }  x_{\pi} c_\pi$ and $\theta_{\gb,m}
=\sum_{\gs }  x_{\gs} \a_\gs$ are related by
\by 
\label{1789} \sum_{\pi\in {{\bf P}}(m\gc)}  (1\ot P_\ga) x_{\pi}\gk( c_\pi)
&=& \sum_{\gs \in {{\bf P}}(m\gb)} ((B_S\ci\nabla_S^t)\ot 1)x_{\gs} \gk( \a_\gs). 
\ey 
\et \bpf 
Apply the left side  of Lemma \ref{dra} (d) to the left side of the  noncommutative version, Equation   \eqref{17x}  and apply the right side  to the right side of 
 \eqref{17x}.
\epf
\noi 
By \eqref{zra1}
 and \eqref{zra2} we have an isomorphism 
\be \label{2W4}f: R_\Z(\gb,\gb)
\lra R_\Z(\gc,\gc).\ee where
$$f:=(1\ot P_\ga^{-1})(B_S\ci\nabla_S^t\ot 1) = (B_S\ci\nabla_S^t\ot 1)
(1\ot P_\ga^{-1}).
\ff{Here $(1\ot P_\ga^{-1})=(1\ot P_\ga)^{-1}$ is the inverse of the map defined in \eqref{zra1}.}$$
Then by Theorem \ref{zra} 
\be \label{za1} 
f( \gk(\theta_{\gb,m})) =\gk(\theta_{\gamma,m}).\ee
Taking into account the gymnastics necessary to define $B_S$, 
we note that $f$ factors through $F$, where  
\be \label{2W9}F=(1\ot P_\ga^{-1})(B_S\ci\nabla_S^t\ot 1):
R_\Z(\gb,\gb)\lra S_\Z(\fn^-)_e \ot \cO_\Z(H_{\gc,0}).\ee
By induction and \eqref{za1}, we have 
\bc  \label{3W} \bi \itema If $\gc$ is isotropic then  $\theta_{\gamma}^2 =0.$
\itemb If $\gc$ is nonisotropic and $m$  satisfies condition (b) before  Theorem \ref{1sag} then $\gk(\theta_{\gamma,m}) = \gk(\theta_{\gamma,1})^m$.
\itemc For any rogue partition 
$\pi$ of $m\gc$ the coefficient of $e_{-\pi}$ in $\theta_{\gc,m}$ is zero.
\ei
\ec \bpf Since $f$ is an algebra map, this follows since $\gk( c_\pi) =0$  implies  $ c_\pi =0$. \epf
\noi See \cite{Md} for a different approach.
\bc \label{3W2} We have equivalences of categories and a commutative diagram. 
$$
\xymatrix{{\bf gr} R(\gb,\gb)/\tors   \ar@/_1pc/[rr]_{{(B_S\ci\nabla_S^t\ot 1)^{*}}} 
\ar@{>}^{f^{*}}[r]
&{\bf gr} R(\gc,\gc)/\tors  \ar@{>}^{{(1\ot P_\ga^{-1})}^{*}}[r]&{\bf gr}R(\gc,\gb)/\tors .}$$ 
\ec

\subsection{Finite generation}  
\label{2W} \noi The next result might appear redundant since we mainly work with small Šapovalov algebras, which are generated by their  space of elements of degree $\le 1$. 
However we could expect to obtain other results by applying convex geometry to
Šapovalov algebras.  For example they should be Cohen-Macaulay rings, \cite{BH} Chapter 6.
Assume that $\fg $ is finite dimensional. 
\bl \label{2W1} The algebras Š$_\gc$  and Š$_\gb$  are finitely generated. 
\el
\bpf We prove the result for Š$_\gc$.  The first step is to show that Š$_\gc$ is isomorphic to a certain semigroup algebra $\ttk S$.  Let 
$S(m) =\{(k,\pi)\in {\Z\ti
\bf{\overline{P}}}(m\gc-k\ga)\},$ using the notation 
from Lemma \ref{1uebasis}.
Then $S= \op_{m \ge} S(m)$ is an additive semigroup, such that Š$_\gc=\ttk S$ the semigroup algebra.
 Let $\cA=\{\ga_1,\ldots \ga_R\}$ be the set of positive roots such that for each $i$ there exists an integer $k$ we have $\ga_i < k\gc$. This means that there is a partition of  $k\gc$ involving $\ga_i.$
Let $M=\sum_{i=1}^R \Z v_i$ be a free abelian group of rank $R$. 
 Assume $\ga=\ga_R$ and set 
$$H= \R v_R +\sum_{i=1}^{R-1} \R_{\ge 0} v_i \subset \R^R =M\ot \R.$$
For $t\in  \R_{\ge 0}$,  set 
$$F_{t}=\{\sum a_i v_i\in H| \sum a_i \ga_i =t\gc\}.$$ 
and  consider the convex set  $F:=\bcu_{t\in  \R_{\ge 0}} F_t.$  
For $m\in \Z^+$ and  a pair $(k,\pi)\in {\Z\ti
\bf{\overline{P}}}(m\gc-k\ga)$ we associate the point 
$v_{k,\pi} = k v_R + \sum_{i=1}^{R-1} \pi(\ga_i) v_i \in F\cap M.$ 
Since each element of $F_m\cap M$  corresponds to a partition of $m\gc,$
this gives a surjective map from the basis of  $S_\gc(m)$ given in Lemma \ref{1uebasis} to $F_m\cap M$ which extends to an homomorphism of semigroups $\phi:S\lra F\cap M.$  \\ 
\\
We show $F\cap M$ is finitely generated. For this we  adapt the proof of Gordan's Lemma \index{Gordan's Lemma} from \cite{Ful} Chap 1. Prop. 1.
 The  set $K=\bcu_{m>0} \frac{1}{m}F_m \cap M $ is a discrete subset of the compact Hausdorff space $F_1$, so is finite. 
Therefore for some $N$, we have $NK\subseteq (F_N \cap M )$.  For example we can take $N$ to be the least common multiple of a finite set of denominators that are required to obtain $K$. Note that $NcK \subseteq M$ for all $c\ge 0.$
We claim that $F\cap M $ is generated by the finite dimensional subspace $\sum_{m \ge 0}^N F_m\cap M .$
Indeed, 
if $v\in F_a\cap M$, then $v/a = s \in K\subset F_1$. Write $a =N c +r$ where $c\ge0$ and  $0\le r < N$. Then
$v = N cs  +rs $. Since $v, N cs \in M$, it follows that $rs\in F_r\cap M$. This proves the claim.\\ \\
The map $\phi$ extends to a map of free abelian groups 
$\phi:\Z S\lra \Z(F\cap M)$. 
Let $B=\Ker \phi 
$ 
This is a free abelian group which we identify $B$ with $\Z^Y$ for some $Y$. Then for $I\subseteq Y$, let $$B_I=
\{(b_1 \ldots, b_Y)| b_i \ge 0 \mbox{ iff } i\in I\}.$$ Then $B=\bcu B_I$.
Each $B_I$ is a finitely generated semigroup (Gordan's Lemma again).
Now let $s_1, \ldots ,s_Z$ be a finite set of elements of $S$ such that 
$F\cap M$ is generated by the image of this set and $\{s_1, \ldots ,s_Z\}$ contains a set of generators for each $B_I$. We claim $s_1, \ldots ,s_Z$ generates $S$. If $x\in S$, then $\phi(x) = \sum_i r_i \phi(s_i)$ where $r_i$ is a nonnegative integer. Since  $x-  \sum_i r_i s_i \in B_I$ for some $I$ the result follows easily.  
\epf
 \section{Geometry of Šapovalov elements}  \label{GF} 
We have   been doing algebraic geometry all along.
 \ff{ll y a plus de quarante ans que je dis de la prose sans que j'en susse rien, Molière, Le Bourgeois Gentilhomme}
\\ \\
\noi Here we work over $\ttk$ since some of the results we need have a Noetherian hypothesis. 
\subsection{Sheaves and morphisms}  \label{GF1}
Let $\phi:X\lra Y$ be a morphism of schemes and $\cF, \cG$ quasi-coherent sheaves on $X, Y$. The {\it direct image } \index{Direct image sheaf} $\phi_* \cF$ of $\cF$ and the {\it inverse image }  \index{Inverse image sheaf}
$\phi^* \cG$ of $\cG$ are defined in \cite{Ha} page 109. We have a natural isomorphism 
 \be \label{XL}\Hom_{\cO_X}(\phi^* \cG, \cF) 
\cong  \Hom_{\cO_Y}(\cG, \phi_*\cF).\ee  
The category of affine schemes is the opposite category to the category of commutative rings. If $f¨:R\lra A$ is a ring homomorphism,  let $X=\Spec A$, $Y=\Spec R$. Then we have a corresponding  morphism of schemes $\phi:X\lra Y$.
 For $N\in R$-mod and $M\in A$-mod, we have corresponding sheaves $\tilde N$, $\tilde M$ on $Y, X$ respectively,
 and \eqref{XL} becomes. 
\be \label{XLL}\Hom_{\cO_X}(\phi^* \tilde N, \tilde M) 
\cong  \Hom_{\cO_Y}(\tilde N, \phi_*\tilde M).\ee
Thus \eqref{XL} is the analog for sheaves of \eqref{GF11}. 

\subsection{A theorem of Serre}

 Given an $\N$ graded ring $R=\bop_{n\ge 0} R(n)$ with $R(0)$  a finitely  generated  $\ttk$-algebra,  such that $R$ is generated over  $R(0)$ by $R(1)$,
 we can construct a scheme $X=\pro R$.  The construction gives invertible sheaves $\cO_X(n)$ on $X$, \cite{Ha} II.2 and II.5.   
Denote the category of quasi-coherent sheaves on $X$ by $\frak{Qco}X$.
\index{$\frak{Qco}X$, Category of quasi-coherent sheaves on $X$}     Then by a Theorem of Serre \cite{Se}, $\frak{Qco}X$ is equivalent to the quotient category ${\bf gr} R/\tors$. 
\\\\
In   \eqref{zr6} and \eqref{zra4}  we introduced 
graded rings 
$R(\gt,\gs)$ 
and various isomorphisms between them, see also \eqref{zra1}, \eqref{zra2}. Set $\X_{\gt,\gs}= \pro R(\gt,\gs)$ \index{Algebras $R(\gt,\gs)$!$\X_{\gt,\gs}= \pro R(\gt,\gs)$}. To avoid an accumulation of notation we sometimes use the same symbol for a map between graded rings and the corresponding map between projective schemes (when it exists).

\bp
\bi \itema We have a commutative diagram of projective schemes

\xymatrix{
&&&&\X_{\gc, \gb}
 \ar@{->}[dr]_
{(1\ot P_\ga^{-1})}
  \ar@{->}[rr]^{(B_S\ci\nabla_S^t\ot 1)} 
&& \X_{\gb, \gb}
\ar@{<-}^{f
}[dl] &\\
&&&&&\X_{\gc, \gc}
} 
\itemb
\noi We have equivalences of categories and a commutative diagram. 
$$
\xymatrix{\Qbb  \ar@/_1pc/[rr]_{{(B_S\ci\nabla_S^t\ot 1)^{*}}} 
\ar@{>}^{f^{*}}[r]
& \Qcc \ar@{>}^{{(1\ot P_\ga^{-1})}^{*}}[r]&\Qcb}.$$ \ei
\ep\bpf (a) follows since the isomorphisms preserve the irrelevant ideals and (b) holds by 
 Corollary  
\ref{3W2} and the Theorem  of Serre.\epf
\noi By zooming in, the map $f$ from (a) induces an isomorphism of fiber products over $\A^1.$
$$\X_\gc \ti \A^1\lra \X_\gb \ti \A^1.$$
Let $\cO_{\gt,\gs}(m) =\cO_{\X_{\gt,\gs}}(m).$ Then sections of 
$\Obb(m)$ have the form
$\sum_{\gs \in {{\bf P}}(m\gb)}x_{\gs} \gk(\a_\gs).$
\\ \\Let $Z = (\fn^-)^*$, so that $S(\fn^-)=\cO(Z) $ and $S(\fn^-)_e=\cO(Z)_e = \cO(Z_e)$ where 

\be\label{GF7} 
Z_e = \{z\in Z|e(z)\neq 0\}
\ee
Note that 
$\Obb(m)$ and  $(\nabla_S)_*\Obb(m)$  are 
objects in $\Qbb$ and  $\Ocb(m)$ 
is an object in 
$\Qcb$ and we have maps 
\be\label{318}((B_S\ci\nabla_S^t)\ot 1)^{*}:
 \Obb(m){\lra}  \Obb(m) \lra \Ocb(m)\ee
\be\label{319}(1\ot P_\ga^{-1})^{*}:
\Occ(m)  {\lra}  \Ocb(m)\ee
 Now  $$\gk( \theta_{\gamma,m})=\sum_{\pi \in {{\bf P}}(m\gamma)} 
x_{\pi}\ot \gk( c_\pi) 
\mbox { and }  
\gk( \theta_{\gb,m})=\sum_{\gs \in {{\bf P}}(m\gb)} x_{\gs}\ot \gk(\a_\gs) 
$$
are sections of $\Occ(m), \Obb(m)$ respectively. 
If we evaluate \eqref{318}
and \eqref{319} on $ 
\gk( \theta_{\gb,m})$ and $\gk(\theta_{\gamma,m})$ respectively,
we find exactly the same relation \eqref{1789} as in
Theorem  \ref{zra}, except that algebra maps have been replaced by
inverse image functors. Moreover by \eqref{za1},
\be \label{zb1} 
f^*( \gk(\theta_{\gb,m})) =\gk(\theta_{\gamma,m}).\ee

\subsection{The cancellation step}
By \eqref{2W9}, there is a map $F:R(\gb,\gb)\lra S(\fn^-)_e \ot \cO(H_{\gc,0}). 
$ The cancellation step in Subsection   \ref{cst}  shows 
$F(\gk(\theta_{\gb, m}))(\gl)\in S(\fn^-)$ for all   
$\gl\in H_{\gc,m}$. If $\gth\in\cO(Y)(m)= R(\gb,\gb)(m)$, 
the statement $F(\theta)(\gl)\in S(\fn^-)$ can be interpreted as follows.
\\ \\
Recall the map 
$y_m^*:H_{\gs, 0} \lra H_{\gs, m}$ from \eqref{jjj}
\bp 
Given $\gl  \in H_{\gc,0}$ 
\bi \itema There is a graded subalgebra 
$A_\gl= \bop_{m\ge 0}A_\gl(m)$ of $\mbox{Š}^0_{\gb}$ with
$$A_\gl(m)\ot\cO(H_{\gb,m})=\{\gth\in R(\gb,\gb)(m)|(F\gth)(y_m^*\gl)\in S(\fn^-) \}.$$
\itemb 
$A_\gl(m)$ contains $\gth_{\gb,m}$.
\ei\ep \bpf Use Lemma \ref{td5} for (a) and the cancellation step for (b). \epf 
\subsection{Existence and uniqueness of \v Sapovalov elements.} \label{inc}
\bl Given $\gth\in \Obb(m)$ and $S \subseteq \Gd^+$
$$C_{\gth, S} =\{\mu\in \fh^* |e_\ga \gth v_\mu =0 \mbox{ for all } \ga\in S\}$$
 is closed in $\fh^*$.
\el \bpf This holds by the proof of \cite{D} Lemma 7.6.12.\epf
\noi If $S = \Gd^+$  and $\theta=\theta_{\gb,m}$, then $C_{\gth, S}= {H}_{\gb, m}$\bp 
 Given $S$ and $C$ a closed set in $\fh^* $
define $$W_{S, C} =\{\gth\in \Obb(m)|e_\ga \gth v_\mu =0 \mbox{ for all } \ga\in S\}$$
If $S= \Gd^+$ and $C = {H}_{\gb, m}$  then $W_{S, C} $ is one dimensional and spanned by $\theta_{\gb,m}$.\ep  
\bpf This holds by the existence and uniqueness of \v Sapovalov elements.
\epf
\section{{Examples}}  \label{vf}  
\subsection{Type A}  \label{Type A}  
\bexa \label{PA.13}{\rm
This is one of the smallest examples of a map between \v Sapovalov  algebras $$R=\mbox{Š}_\gb \lra  \mbox{Š}_\gc =S.$$  It arises in the computation of the \v Sapovalov element $\theta_{\gamma,1}$ for the highest root $\gc$ of $\fg=\fsl(3)$.  
Let $\fg' \cong \fsl(2)$ be the subalgebra of $\fg$ consisting of $2\ti 2$ matrices located in the lower right corner of $\fg.$   Let $\gb$ be the simple root of $\fg'$ and $\ga$ the other simple root of $\fg$. Then 
$ X_\ga =Y, X_\gb =X, Z=[X,Y]=X_\gc. $
form a basis for $\fn^-$.
\\  \\ The \v Sapovalov subalgebras are 
 $R =\ttk[X, ZY^{-1}], S=\ttk[XY, Z],$ with 
$$\deg_R X =\deg_R ZY^{-1}   =1,  \deg_S XY =\deg_S Z = 1 \mbox{ and } \deg T =0.$$ 
We have $d_- = Y^{-1}\ad Y$, and 
$\nabla_S = \exp(Td_-)=1  +Td_- + \mbox{ terms with higher powers of }Td_- $. Since $d_-(X) = -ZY^{-1}$, $\nabla_S$ is the automorphism of $R[T]$ given by 
$$ X\lra X -TZY^{-1}, \quad ZY^{-1}\lra ZY^{-1}.$$
\noi
The composite  $B \ci \nabla_S: R[T]\lra S[T]$ satisfies 
$$ X\lra XY -TZ, \quad ZY^{-1}\lra Z.$$
To compute the 
  \v Sapovalov element $\theta_{\gamma,1}$ we use Equation 
\eqref{12nd}.  This yields 
$$\theta_{\gamma,1}(s_\ga\cdot\mu)= e_{- \gb}e_{- \ga} - te_{- \gc}.$$}\eexa
\subsection{\v Sapovalov elements for $\fsp(4)$} \label{g19}
Let $\fg =\fsp(4)$ with positive roots $\gb,  \ga+\gb, 2\ga +\gb, \ga$.  Denote the corresponding negative root vectors by $e_{-\gb}, e_{-\ga-\gb}, e_{-2\ga -\gb}, e_{-\ga}$. 
\subsubsection{Case 1}  \label{71a} We find the element $\gth_{\ga+\gb,1}$ for $s_\ga(\gb)= \ga+\gb.$ 
For this we invert $e_{-\ga}$ and then look at the $\gb$ weight space. 
There are two basis elements $e_{-\gb}, e_{-\ga-\gb}e_{-\ga}^{-1}$. We compute 
$$\gth_{\ga+\gb,1} = e_{-\ga}e_{-\gb} + (\mbox{ coeff of deg 1})e_{-\ga-\gb}.$$ Thus no rogue partitions occur in $\gth_{\ga+\gb,2} $.
But on the other hand we have an isomorphism of graded algebras $ \mbox{Š}_\gb\lra  \mbox{Š}_{\ga+\gb}$ and $e_{-2\ga -\gb}e_{-\gb}\in \mbox{Š}_{\ga+\gb}(2)$ corresponds to a   rogue partition of $2(\ga+\gb)$.  The explanation for this is that  $ e_{-2\ga -\gb}e_{-\gb}e_{-\ga}^{-2
}$ is an element in degree two  which $B_S$ maps to $e_{-2\ga -\gb}e_{-\gb}$.
\subsubsection{Case 2}  \label{71b}
To find the element $\gth_{2\ga+\gb,1}$ for $s_\gb(\ga)= 2\ga+\gb,$we invert $e_{-\gb}$ and look at the $\ga$ weight space.
There are three basis elements $e_{-\ga},\; e_{-\ga-\gb}e_{-\gb}^{-1},\; e_{2\ga +\gb}e_{-\gb}^{-2}$. We compute 
$$\gth_{2\ga+\gb,1} = e_{-\ga}e_{-\gb}^2 + (\mbox{ coeff of deg 1})e_{-\ga-\gb}e_{-\gb} +(\mbox{ coeff of deg 2})e_{-2\ga -\gb}.$$ 
\subsection{Type C: The cancellation step}  \label{Type C}
{\rm 
\noi A  crucial step in the construction of \v Sapovalov elements is the observation
that the terms 
$ q^{(j)} _{\gs}( e_{-\gz}^\ttk)$ in Lemma 
\ref{109b}  
vanish unless 
\begin{eqnarray} \label{198g}
\gz(\ga)\ge N- mq.\nn
\end{eqnarray}
 This follows from representation theory. The condition is necessary since otherwise the evaluation of the \v Sapovalov element would belong to 
$U(\fn^-)_e$ and not to  $U(\fn^-)$.
We give an example where the individual terms $R _{\gs, j}$ on the right of 
Equation (\ref{a1}) are not identically zero, and verify directly that  the sum in Lemma \ref{1c} (b) is zero, see \eqref{vanish}. 
Consider the Dynkin-Kac diagram below for the Lie superalgebra $\fg = \osp(2,4)$.

\vspace{0.4cm}

\vspace{0.4cm}
\begingroup
\setlength{\unitlength}{0.10in}
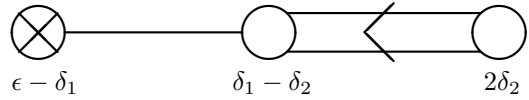
\begin{picture}(-30,-10)
\thicklines
\put(14.414,0.0){\line(1,0){9.23}}
\put(25.95,1.0){\line(1,0){10.1}}
\put(25.95,-1.0){\line(1,0){10.1}}
\put(30,0){\line(1,-1){1.5}}
\put(30,0){\line(1,1){1.5}}
\put(12,1){\line(1,-1){2}}%
\put(12,-1){\line(1,1){2}}
\put(13,0){\circle{2.828}}
\put(25,0){\circle{2.828}}
\put(37,0){\circle{2.828}}
\put(11.6,-3.0){$\gep-\gd_1$}
\put(23.1,-3.0){$\gd_1-\gd_2$}
\put(36.2,-3.0){$2\gd_2$}
\end{picture}

\vspace{1.4cm}

\endgroup
\noi If we replace the grey node in the diagram  labelled by the root $\gep-\gd_1$ by a white node, labelled by the root $\ga_0 =\gd_0-\gd_1$ we obtain the Dynkin diagram and root system of $\fsp(6)$.
Note that the equations below involving $\gr$,  \eqref{k1} - \eqref{k3} do not depend on 
$\gep-\gd_1$ or $\ga_0 =\gd_0-\gd_1$. Although the dimensions of the nilradicals 
have dimensions 8 and 9 for $\fg = \osp(2,4)$ and   $\fsp(6)$ respectively, we do not need the root $2\gd_0$ in any  of the computations below.
Another point to make is that  for the case $\fg = \osp(2,4)$ in any partition $\pi$ that we use, there is at most one odd root  $\gc$ such that $\pi(\gc)> 0$. Therefore the commutation relations needed for 
 $\fg = \osp(2,4)$ are exactly the same as  $\fsp(6)$.\\ \\
\noi Let $\fb$ be the Borel subalgebra having  simple
roots  $\gb = \epsilon - \gd_1, \ga_1 = \gd_1 - \gd_2,$ $  \ga_2 = 2\gd_2$ and 
$\fb^-$  the Borel opposite to $\fb$. In the case of 
$\fsp(6)$ set  $\gb = \gd_0- \gd_1.$ 
 Let
$e_{-\gb}, e_{-\ga_1}, e_{-\ga_2}  $ be the negative simple root vectors.
We compute the \v Sapovalov elements
$\gth_1, \gth_2, \gth_3$ for the positive roots $\gc_1= \gb + \ga_1, \gc_2= \gb + \ga_1+ \ga_2$ and $\gc_3= \gb + 2\ga_1+ \ga_2$.
\subsubsection{\v Sapovalov elements}  \label{xy.2}
Define the other negative root vectors by
$$e_{- \ga_1 - \ga_2} = [e_{-\ga_1},e_{- \ga_2 }  ],  \quad \quad
e_{- 2\ga_1 - \ga_2} = [e_{-\ga_1}, e_{- \ga_1 - \ga_2}],  $$
$$e_{-\gb- \ga_1} =[e_{-\ga_1} ,e_{-\gb}
], \quad e_{-\gb- \ga_1-\ga_2} = [e_{- \ga_2 } ,e_{-\gb- \ga_1}], \quad e_{-\gb- 2\ga_1 - \ga_2} = [e_{-\ga_1},e_{-\gb- \ga_1-\ga_2}]. $$
It follows from the Jacobi identity that
$$[e_{-\gb}, e_{- \ga_1 - \ga_2}] = e_{-\gb- \ga_1-\ga_2},
\quad \quad
[e_{- \ga_1 - \ga_2},e_{-\gb- \ga_1}] = e_{-\gb- 2\ga_1 - \ga_2},$$ and
$$[e_{-\gb},e_{- 2\ga_1 - \ga_2}] = 2e_{-\gb- 2\ga_1-\ga_2}.$$
\noi 
Throughout we order the set of positive roots so that 
if $\gs$ belongs to the set
\be\label{X1}\{
\gb, \gb + \ga_1, \gb + \ga_1 + \ga_2, \gb + 2\ga_1 + \ga_2\}\ee 
if $\pi$ is a partition   
with $\pi(\gs)> 0$ then  $e_{-\gs}$ is the 
first factor in $e_{-\pi}$. There is a unique $\ga_1$ order and a unique $\ga_2$ order satisfying these conditions. In \eqref{1la2} we switch the orders.
\\ \\
Let $s_1, s_2$ be the reflections corresponding to the simple roots $\ga_1, \ga_2.$
Then for $\gl \in \fh^*$ define $\gl_1 =s_1\cdot\gl, \;\gl_2 = s_2\cdot\gl_1, \;\gl_3 = s_1\cdot\gl_2$. Using the notation \eqref{efh}, set $t_i= t_{\ga_i}$ for $i=1,2$. For a root $\gc$ any polynomial function on ${H}_{\gc}$  can be expressed uniquely as a polynomial in $t_1, t_2.$
Let
\be\label{k1}p=(\gl + \gr,\ga^\vee_1) = 
-(\gl_1 + \gr,\ga^\vee_1) =  -(\gl_2 + \gr,(\ga_1+\ga_2)^\vee)=-(\gl_3 + \gr,(\ga_1+\ga_2)^\vee),\ee 
\be\label{k2}
q= (\gl + \gr,(2\ga_1  + \ga_2)^\vee) = (\gl_1 + \gr,\ga^\vee_2) = -(\gl_2 + \gr,\ga_2^\vee) = -(\gl_3 + \gr,(2\ga_1+\ga_2)^\vee) ,\ee 
and 
\be\label{k3}r=(\gl + \gr,(\ga_1+\ga_2)^\vee) =(\gl_1 + \gr,(\ga_1+\ga_2)^\vee)  =
(\gl_2 + \gr,\ga^\vee_1)  = -(\gl_3 + \gr,\ga_1^\vee) .\ee
Then $r = 2q - p.$  
We compute the
\v Sapovalov elements
$\gth_{\gc_i}$ for $\fg$. 
To do this we use Equation (\ref{121nd}) and  the abbreviations $\gth_{i}: = \gth_{\gc_i}(\gl_i)$.   
To start the induction, assume that $(\gl+\gr,\gb) =(\gb,\gb)/2.$ Suppose that $p, q, r$ are nonnegative integers. Then by \eqref{121nd}
\be \label{nelab}
e_{-\ga_1}^{p+1}e_{-\gb} = \gth_1e_{-\ga_1}^{p},\quad e_{- \ga_2 }^{q+1}\gth_1 = \gth_2e_{- \ga_2 }^{q},\quad 
e_{- \ga_1 }^{r+1}\gth_2 = \gth_3e_{- \ga_1 }^{r}.\ee 
In the Verma module $M(\gl)$ induced from $U(\fb)$ with highest weight $\gl$ and highest weight vector $v_\gl$, set $v_1=e_{-\ga_1}^{p}v_\gl,$ $v_{ _2  }=e_{- \ga_2 }^{q}v_1$ and $v_{3}= e_{- \ga_1 }^{r}v_2$. These are all highest weight vectors for 
$\fb$ and by \cite{M101}, Theorem 9.3.2 they generate Verma submodules of $M(\gl)$. 
In $U(\fb^-)$ we have 
$$[e_{-\ga_1}^{p+1} ,e_{-\gb}
] = (p+1)e_{-\gb- \ga_1} e^p_{-\ga_1}$$
$$[e_{- \ga_2 }^{q+1} ,e_{-\gb- \ga_1}] =
(q+1)e_{-\gb- \ga_1-\ga_2} e_{- \ga_2 }^q$$
$$[e_{- \ga_2 }^{q+1} ,e_{-\ga_1}]  =-(q+1)e_{- \ga_1 - \ga_2}
e_{- \ga_2 }^q.$$
This easily gives
\be \label{1tha}
\gth_1 = (p+1)e_{-\gb- \ga_1} + e_{-\gb}e_{-\ga_1} = pe_{-\gb- \ga_1} + e_{-\ga_1}e_{-\gb}.\ee
We order the set of positive roots so that for any  partition $\pi$, $e_{-\ga_2}$ occurs last if at all in $e_{-\pi}$. Then %
\begin{eqnarray}\label{1la2}
\gth_2
&=& (p+1)[(q+1)e_{- \gb - \ga_1 - \ga_2} +e_{-\gb- \ga_1}e_{-\ga_2}] +e_{-\gb}[ e_{-\ga_1}e_{- \ga_2} - (q+1)e_{- \ga_1 - \ga_2} ].\nn\\
&=& (p+1)[(q+1)e_{- \gb - \ga_1 - \ga_2} +e_{-\gb- \ga_1}e_{-\ga_2}] +e_{-\gb}[ e_{-\ga_2}e_{- \ga_1} - qe_{- \ga_1 - \ga_2} ].
\end{eqnarray}
Next order the set of positive roots so that for any  partition $\pi$, $e_{-\ga_1}$ occurs last $e_{-\pi}.$ To find $\gth_3$  we use
$$[e_{-\ga_1}^{r+1} ,e_{-\gb-\ga_1-\ga_2}
] = (r+1)e_{-\gb- 2\ga_1 -\ga_2} e^r_{-\ga_1},$$
$$[e_{-\ga_1}^{r+1} ,e_{-\gb-\ga_1}e_{-\ga_2}
] = (r+1)e_{-\gb- \ga_1}e_{-\ga_1-\ga_2}e^
{r}_{-\ga_1}
+\left( \begin{array}{c}
                r+1 \\
                2 \end{array}\right)e_{-\gb- \ga_1}e_{-2\ga_1 -\ga_2} e^{r-1}_{-\ga_1},$$
$$[e_{-\ga_1}^{r+1} ,e_{-\gb}e_{-\ga_1-\ga_2}
] = (r+1)[e_{-\gb}e_{- 2\ga_1 -\ga_2} e^r_{-\ga_1}
 +e_{-\gb- \ga_1}e_{-\ga_1 -\ga_2} e^r_{-\ga_1} +re_{-\gb-\ga_1}e_{-2\ga_1-\ga_2}
e_{\ga_1}^{r-1}],$$
$$e_{-\ga_1}^{r+1} e_{-\gb}e_{-\ga_2}e_{-\ga_1}
 =
e_{-\gb}[e_{-\ga_2}
e_{- \ga_1}^{2} +
(r+1)e_{- \ga_1-\ga_2}e_{-\ga_1}
+\left( \begin{array}{c}
                r+1 \\
                2 \end{array}\right)e_{- 2\ga_1 -\ga_2}] e^{r}_{-\ga_1}$$
$$+ (r+1)e_{-\gb- \ga_1}[e_{-\ga_2}
e_{- \ga_1} +
re_{- \ga_1-\ga_2}]e^{r}_{-\ga_1}
+(r-1)\left( \begin{array}{c}
                r+1 \\
                2 \end{array}\right)e_{-\gb-\ga_1}e_{-2\ga_1-\ga_2}
e_{\ga_1}^{r-1}
                .$$
The above equations allow us to write $e_{- \ga_1 }^{r+1}\gth_2$  in terms of elements $e_{-\pi}$ with $\pi$ a partition of $\gb + (r+2)\ga_1 + \ga_2$.
We see that the term $e_{-\gb- \ga_1}e_{- 2\ga_1 -\ga_2} e^{r-1}_{-\ga_1}$  occurs in $e_{- \ga_1 }^{r+1}\gth_2$ with coefficient
\be \label{vanish}
\left( \begin{array}{c}
                r+1 \\
                2 \end{array}\right)[(p+1) -2q + (r-1)] = 0.\ee
This is predicted by the cancellation step in Subsection \ref{cst}. 
In the remaining terms, $e_{- \ga_1 }^{r}$ can be factored on the right, and this yields
\begin{eqnarray}\label{1th3}
\gth_3 &=&
(p+1)(q+1)(r+1)e_{-\gb- 2\ga_1 -\ga_2} +(p+1)(q+1)e_{-\gb- \ga_1-\ga_2}e_{-\ga_1} \\
&+& (q+1)(r+1)e_{-\gb- \ga_1}e_{-\ga_1 -\ga_2}
-(p/2)(r+1)e_{-\gb}e_{- 2\ga_1 -\ga_2}\nonumber \\
&+& 2(q+1)e_{-\gb- \ga_1}e_{-\ga_2}e_{-\ga_1}+ (r-q+1)e_{-\gb}e_{-\ga_1-\ga_2}e_{-\ga_1}+e_{-\gb}e_{-\ga_2}
e_{- \ga_1}^{2}. \nonumber
\end{eqnarray}

\subsection{Type C: Decomposition into simple $\fl$-modules}  \label{Type C2}
Given the importance of the adjoint action of $\fl$, we might hope 
to decompose each  graded component  $\mbox{Š}_{\gb,\Z}(m)$  is a direct sum of simple $\fl$-modules (over $\Z$). 
We focus on the $\fsl(2)$ case and indicate how this can be done, see Example \ref{nB2} for the case of $\fsp(6)$ when $mq=1$. However the process depends on the choice of a subset of 
${{\bf P}}(m\gb)$ and there can also be highest weight vectors which have poles, see Example \ref{nB2} for the case of $\fsp(6)$ when $mq=1$.
Hence outside of Type A see \cite{M23}, it seems unlikely that there is a canonical way to do this.

\subsubsection{How to improve a $\Z $-basis}   Let $\fsl(2)$ be the  Lie algebra with basis $E, F, H$ such that $[H,E]=2E, [H,F]=-2F$ and  $[E,F]=H$. Denote the simple $\fsl(2)$ module of dimension $n+1$ by  $L(n)$. In the adjoint representation, $E, F$ act as $d_{+}, d_{-}$ respectively and $H$ acts as $\ad h_{\ga}.$  
We want to improve the $\Z $-basis  given by Lemma \ref{1uebasis} 

\bt \label{Z-b} Let  \be M=\bigoplus_{n\in \Z} M(n)\ee  be a $\Z$ lattice in a finite dimensional $\fsl(2)$-module and $u\in M(n)$. 
Assume that $M$ is stable under $E^{(k)}, F^{(k)}$ for all $k$. 
Then there is a unique expression 
\be \label{arl}u = \sum_{k\ge 0} F^{(k)} u(k)\ee where $u(k) \in M(n+2k)\cap \ker E.$
\et 
\bpf This is the limit as $q\lra 1$ of \cite{HK} Lemma 4.1.1.  It is not hard to give a direct proof.
\epf 
\noi The Theorem  shows  
any weight vector $v$ of positive weight is congruent to a highest weight vector  $v_0$ mod the image of $F$. 
\noi For $\cB \subseteq M$, let $\langle \cB \rangle$ be the $\fsl(2)$ submodule generated by $\cB$-

\bc \label{ZB} For all $k>0$ 
 suppose $\cB_k$  is a $\Z$-basis for $M(n+2k)\cap \ker E$ and set $\cB= \ds \cB_k$. Let $\cC_0$ be a $\Z$-basis for $M(n)$. 
Then there is a $\Z$-basis $\cB_0$ of $M(n)$ consisting of  highest weight vectors such that 
$\langle \cB_0 \cup \cB \rangle= \langle \cC_0 \cup \cB \rangle$.
\ec
\bpf This follows directly from Theorem \ref{Z-b}.
  \epf \noi This allows us to write  $M$ as a direct sum of simple modules.  We start by collecting the simple modules of largest dimension and then proceed to lower dimension.  The process is probably best understood from an example. 
\bexa \label{nB2} {\rm 
Let  $\fg =\fsp(6)$ with simple roots $\ga_0, \ga_1, \ga_2$ where $\ga_2$ is a long root and $\ga_0+\ga_2$ is not a root.   Let $x_0, x_1, x_2$ be root vectors with $\wt\, x_i = -\ga_i$.  If $\gb=\ga_0+ \ga_1+ \ga_2$ and $\gc=\ga_0+ 2\ga_1+ \ga_2$, we decompose $M=$Š$_\gb(1)$ into a direct sum of simple modules. This example is essentially the same as in Subsection \ref{Type C}  where $\gb, \gc$ were called $\gc_2$ and $\gc_3$. 
 We have $mq=1.$
 There are 4 partitions $\gs$ of $\gb$ and we order the corresponding $x_\gs$ as    
\be \label{1bc} x_0 x_1x_2,\quad [x_0, x_1]x_2,\quad x_0 [x_1,x_2],\quad [x_0, [x_1,x_2]].\ee 
 If $x_\gs$ is the first of these, then $x_\gs$
 is a highest weight vector of weight 3 since $d_+(e_{-\ga}) =0$.
 From this   
we obtain the summand $\langle x_0 x_1x_2 \rangle$
 which is isomorphic to $L(3)$.  The last 
 $w=[x_0 [x_1,x_2]]$  
is a highest  vector so  $w_0=w$ and 
$M$ contains $\langle  w_0 \rangle.$ 
Consider the second and third terms  $u=[x_0, x_1]x_2,\; v=x_0 [x_1,x_2]$.   In the notation of \eqref{arl}, $u_0=\pm v_0$ 
and we need to choose one of $u_0, v_0$ to include in the sum.  Let us choose $u_0$. Then 
\by \label{107bc} \langle x_0 x_1x_2 \rangle  \op \langle u_0  
\rangle \op \langle w_0  \rangle \subset M. \ey 
We show that $M$ contains one more submodule of dimension 2. This is most easily done by looking for lowest weight vectors. By Lemma \ref{1uebasis}, $(S(\fn^-)_e)^{-\gb}$ has a basis consisting of all products 
$e_{-\pi} e^k$ 
where  $k \in \mathbb{Z}$ and $\pi \in {\overline{\bf P}}(\gb -k\ga)$ with  $\pi(\ga) = 0.$  It is easily seen that ${\overline{\bf P}}(\gb -k\ga)$ is empty  unless $k\in [-2,1]$.
If $k\ge 0,$ then there are 4 partitions with this property and the corresponding monomials are listed in \eqref{1bc}. It remains to consider partitions with poles. 
 If $k\in [-2,-1]$ we find 4 partitions 
associated to the monomials below, all with weight $-1$.
$$[x_0, x_1][x_1,x_2]e^{-1},\quad x_0 [x_1, [x_1,x_2]]e^{-1},\quad [x_1,[x_2,[x_0, x_1]] ]e^{-1},  [x_0, x_1][x_1 [x_1,x_2]]e^{-2} $$ Since the sum on the left side of \eqref{107bc} has only a 3 dimensional $-1$ weight space, there is a lowest weight vector $z$ of weight $-1$ and we have 
\by \label{107bd} 
M &=& \langle x_0 x_1x_2 \rangle  \op \langle u_0  
\rangle \op \langle w_0  \rangle \op \langle z  \rangle
\nn\\ & \cong & L(3)\op 3L(1).\nn
\ey  }\eexa

\printindex


\begin{bibdiv}
\begin{biblist}
\label{bib}

\bib{BGG1}{article}{
   author={Bernstein, I. N.},
   author={Gel{\cprime}fand, I. M.},
   author={Gel{\cprime}fand, S. I.},
   title={Structure of representations that are generated by vectors of
   higher weight},
   language={Russian},
   journal={Funckcional. Anal. i Prilo\v zen.},
   volume={5},
   date={1971},
   number={1},
   pages={1--9},
   issn={0374-1990},
   review={\MR{0291204 (45 \#298)}},
}


\bib{BGG2}{article}{
   author={Bernstein, I. N.},
   author={Gel{\cprime}fand, I. M.},
   author={Gel{\cprime}fand, S. I.},
   title={Differential operators on the base affine space and a study of
   ${\germ g}$-modules},
   conference={
      title={Lie groups and their representations (Proc. Summer School,
      Bolyai J\'anos Math. Soc., Budapest, 1971)},
   },
   book={
      publisher={Halsted, New York},
   },
   date={1975},
   pages={21--64},
   review={\MR{0578996 (58 \#28285)}},
}

		

\bib{BH}{book}{
   author={Bruns, Winfried},
   author={Herzog, J\"{u}rgen},
   title={Cohen-Macaulay rings},
   series={Cambridge Studies in Advanced Mathematics},
   volume={39},
   publisher={Cambridge University Press, Cambridge},
   date={1993},
   pages={xii+403},
   isbn={0-521-41068-1},
   review={\MR{1251956}},
}
		


\bib{D}{book}{
   author={Dixmier, Jacques},
   title={Enveloping algebras},
   series={Graduate Studies in Mathematics},
   volume={11},
   note={Revised reprint of the 1977 translation},
   publisher={American Mathematical Society},
   place={Providence, RI},
   date={1996},
   pages={xx+379},
   isbn={0-8218-0560-6},
   review={\MR{1393197 (97c:17010)}},
}


\bib{FG}{article}{
   author={Fioresi, R.},
   author={Gavarini, F.},
   title={Chevalley supergroups},
   journal={Mem. Amer. Math. Soc.},
   volume={215},
   date={2012},
   number={1014},
   pages={vi+64},
   issn={0065-9266},
   isbn={978-0-8218-5300-9},
   review={\MR{2918543}},
   doi={10.1090/S0065-9266-2011-00633-7},
}


\bib{Ful}{book}{
   author={Fulton, William},
   title={Introduction to toric varieties},
   series={Annals of Mathematics Studies},
   volume={131},
   note={The William H. Roever Lectures in Geometry},
   publisher={Princeton University Press},
   place={Princeton, NJ},
   date={1993},
   pages={xii+157},
   isbn={0-691-00049-2},
   review={\MR{1234037 (94g:14028)}},
}

\bib{G3}{article}{
   author={Gorelik, Maria},
   title={On the ghost centre of Lie superalgebras},
   journal={Ann. Inst. Fourier (Grenoble)},
   volume={50},
   date={2000},
   number={6},
   pages={1745--1764 (2001)},
   issn={0373-0956},
   review={\MR{1817382 (2002c:17017)}},
}


\bib{Gk}{article}{ author={Gorelik, Maria}, title={The Kac construction of the centre of $U(\germ g)$ for Lie superalgebras}, journal={J. Nonlinear Math. Phys.}, volume={11}, date={2004}, number={3}, pages={325--349}, issn={1402-9251}, review={\MR{2084313 (2005f:17011)}}, }


\bib{Ha}{book}{ author={Hartshorne, Robin}, title={Algebraic geometry}, note={Graduate Texts in Mathematics, No. 52}, publisher={Springer-Verlag}, place={New York}, date={1977}, pages={xvi+496}, isbn={0-387-90244-9}, review={\MR{0463157 (57 \#3116)}}, }

\bib{HK}{book}{
   author={Hong, Jin},
   author={Kang, Seok-Jin},
   title={Introduction to quantum groups and crystal bases},
   series={Graduate Studies in Mathematics},
   volume={42},
   publisher={American Mathematical Society, Providence, RI},
   date={2002},
   pages={xviii+307},
   isbn={0-8218-2874-6},
   review={\MR{1881971}},
   doi={10.1090/gsm/042},
}

\bib{H}{book}{ author={Humphreys, James E.}, title={Introduction to Lie algebras and representation theory}, note={Graduate Texts in Mathematics, Vol. 9}, publisher={Springer-Verlag}, place={New York}, date={1972}, pages={xii+169}, review={\MR{0323842 (48 \#2197)}}, }





\bib{H2}{book}{
   author={Humphreys, James E.},
   title={Representations of semisimple Lie algebras in the BGG category
   $\scr{O}$},
   series={Graduate Studies in Mathematics},
   volume={94},
   publisher={American Mathematical Society},
   place={Providence, RI},
   date={2008},
   pages={xvi+289},
   isbn={978-0-8218-4678-0},
   review={\MR{2428237}},
}


\bib{J1}{book}{ author={Jantzen, Jens Carsten}, title={Moduln mit einem h\"ochsten Gewicht}, language={German}, series={Lecture Notes in Mathematics}, volume={750}, publisher={Springer}, place={Berlin}, date={1979}, pages={ii+195}, isbn={3-540-09558-6}, review={\MR{552943 (81m:17011)}}, }
\bib{Kac1}{article}{ author={Kac, V. G.}, title={Lie
superalgebras}, journal={Advances in Math.}, volume={26},
date={1977}, number={1}, pages={8--96}, issn={0001-8708},
review={\MR{0486011 (58 \#5803)}}, }



\bib{Kac4}{article}{ author={Kac, V. G.}, title={Laplace operators of infinite-dimensional Lie algebras and theta functions}, journal={Proc. Nat. Acad. Sci. U.S.A.}, volume={81}, date={1984}, number={2, Phys. Sci.}, pages={645--647}, issn={0027-8424}, review={\MR{735060 (85j:17025)}}, }

\bib{KK}{article}{ author={Kac, V. G.}, author={Kazhdan, D. A.}, title={Structure of representations with highest weight of infinite-dimensional Lie algebras}, journal={Adv. in Math.}, volume={34}, date={1979}, number={1}, pages={97--108}, issn={0001-8708}, review={\MR{547842 (81d:17004)}}, }		

\bib{KL}{article}{
   author={Kumar, Shrawan},
   author={Letzter, Gail},
   title={Shapovalov determinant for restricted and quantized restricted
   enveloping algebras},
   journal={Pacific J. Math.},
   volume={179},
   date={1997},
   number={1},
   pages={123--161},
   issn={0030-8730},
   review={\MR{1452529}},
   doi={10.2140/pjm.1997.179.123},
}		


\bib{Mq}{book}{
   author={Marquis, Timoth\'{e}e},
   title={An introduction to Kac-Moody groups over fields},
   series={EMS Textbooks in Mathematics},
   publisher={European Mathematical Society (EMS), Z\"{u}rich},
   date={2018},
   pages={xi+331},
   isbn={978-3-03719-187-3},
   review={\MR{3838421}},
   doi={10.4171/187},
}

\bib{Md}{article}{author={Mudrov, Andrei}, title={Factorization of Shapovalov elements},journal={prepint, arXiv:2202.06220},date={2022},}

\bib{M101}{book}{
   author={Musson, Ian M.},
   title={Lie superalgebras and enveloping algebras},
   series={Graduate Studies in Mathematics},
   volume={131},
   publisher={American Mathematical Society},
   place={Providence, RI},
   date={2012},
   pages={xx+488},
   isbn={978-0-8218-6867-6},
   review={\MR{2906817}},
}


\bib{M23}{article}{
   author={Musson, Ian M.},
   title={Explicit expressions for \v{S}apovalov elements in Type A},
   journal={J. Algebra},
   volume={623},
   date={2023},
   pages={358--394},
   issn={0021-8693},
   review={\MR{4557791}},
   doi={10.1016/j.jalgebra.2022.12.033},
}
	

\bib{Sh}{article}{ author={{\v{S}}apovalov, N. N.}, title={A certain bilinear form on the universal enveloping algebra of a complex semisimple Lie algebra}, language={Russian}, journal={Funkcional. Anal. i Prilo\v zen.}, volume={6}, date={1972}, number={4}, pages={65--70}, issn={0374-1990}, review={\MR{0320103 (47 \#8644)}}, }


\bib{S1}{article}{author={Serganova, V.},title={Kac-Moody superalgebras and integrability},conference={title={Developments and trends in infinite-dimensional Lie theory},},book={series={Progr. Math.},volume={288},publisher={Birkh\"auser Boston Inc.},place={Boston, MA},},date={2011},pages={169--218},review={\MR{2743764 (2011m:17056)}}}
		

\bib{Se}{article}{
   author={Serre, Jean-Pierre},
   title={Faisceaux alg\'{e}briques coh\'{e}rents},
   language={French},
   journal={Ann. of Math. (2)},
   volume={61},
   date={1955},
   pages={197--278},
   issn={0003-486X},
   review={\MR{68874}},
   doi={10.2307/1969915},
}
 \end{biblist}
\end{bibdiv}
\end{document}